\documentclass[11pt]{article}
\usepackage[margin=0.95in]{geometry}
\usepackage{epsfig,amssymb,soul}
\usepackage[normalem]{ulem}

\usepackage{hyperref}
\usepackage{graphicx}
\usepackage{tabularx}
\usepackage{adjustbox}
\usepackage{geometry}
\usepackage{diagbox}
\usepackage{pgfplots}
\usepackage{import}
\usepackage{mathrsfs}
\usepackage{standalone}
\usepackage{multirow}
\newcommand{\etal}{\textit{et al}.~}
\pgfplotsset{compat=newest}
 \usetikzlibrary{plotmarks}
  \usetikzlibrary{arrows.meta}
  \usepgfplotslibrary{patchplots}
  \usepackage{grffile}
  \usepackage{amsmath}
\usepackage{epstopdf}
\usepackage{amsfonts,amsmath}
\usepackage{bm}
\usepackage{setspace}
\usepackage{comment}
\usepackage{siunitx}
\usepackage{mathrsfs}
\usepackage{multirow}
\usepackage{subcaption}
\usepackage[hang,flushmargin,bottom]{footmisc}
\usepackage{mathtools}
\usepackage{lipsum}
\usepackage{array}
\usepackage{enumerate}
\usepackage{ifthen}
\usepackage[affil-it]{authblk}
\usepackage{cite}
\usepackage{soul}
\usepackage{leftidx}
\usepackage{relsize}
\usepackage{braket}
\usepackage{float}
\usepackage{chngcntr}
\usepackage{mathrsfs}
\usepackage{bm}
\usepackage[T1]{fontenc}
\usepackage{pxfonts} 

\usepackage{comment}
\usepackage[utf8]{inputenc}

\title{A two-dimensional fractional-order element-free Galerkin method for nonlocal elasticity and complex domain problems}

\author[1]{Shubham Desai}
\author[1]{Malapeta Hemasundara Rao}
\author[1]{Sai Sidhardh\footnote{\textit{Corresponding author: sidhardh@mae.iith.ac.in}}}
\affil[1]{SMART Lab, Department of Mechanical and Aerospace Engineering, IIT Hyderabad, Kandi, Telangana 502285, India}

\date{}

\begin{document}

\maketitle
\begin{abstract}
This study presents a meshfree two-dimensional fractional-order Element-Free Galerkin (2D f-EFG) method as a viable alternative to conventional mesh-based FEM for a numerical solution of (spatial) fractional-order differential equations (FDEs). The previously developed one-dimensional f-EFG solver offers a limited demonstration of the true efficacy of EFG formulations for FDEs, as it is restricted to simple 1D line geometries. In contrast, the 2D f-EFG solver proposed and developed here effectively demonstrates the potential of meshfree approaches for solving FDEs. The proposed solver can handle complex and irregular 2D domains that are challenging for mesh-based methods. As an example, the developed framework is employed to investigate nonlocal elasticity governed by fractional-order constitutive relations in a square and circular plate. Furthermore, the proposed approach mitigates key drawbacks of FEM, including high computational cost, mesh generation, and reduced accuracy in irregular domains. The 2D f-EFG employs 2D Moving Least Squares (MLS) approximants, which are particularly effective in approximating fractional derivatives from nodal values. The 2D f-EFG solver is employed here for the numerical solution of fractional-order linear and nonlinear partial differential equations corresponding to the nonlocal elastic response of a plate. The solver developed here is validated with the benchmark results available in the literature.  While the example chosen here focuses on nonlocal elasticity, the numerical method can be extended for diverse applications of fractional-order derivatives in multiscale modeling, multiphysics coupling, anomalous diffusion, and complex material behavior.\\

\noindent
\textbf{\textit{Keywords: Fractional calculus; Nonlocal elasticity; Element-free Galerkin method; Nonlinear structural mechanics.}}
\end{abstract}
\section{Introduction}
Over the past few decades, fractional calculus has advanced considerably in both theory and applications, motivated by the inability of classical integer-order models to capture anomalous physical behavior. Researchers have proposed novel definitions of fractional derivatives and integrals, supported by rigorous mathematical analysis. Unlike traditional operators, fractional derivatives with their differ-integral definition are intrinsically multiscale, capturing memory effects through temporal components and long-range spatial interactions via nonlocal formulations \cite{oldham1974fractional}. This makes fractional calculus effective for modeling nonlocal and multiscale phenomena across science and engineering. Furthermore, the development of analytical and numerical methods for solving fractional-order differential equations (FDEs) has become an active area of research, enabling accurate simulations of complex systems.

The application of fractional-order spatial and temporal derivatives has expanded considerably within physics, chemistry, and applied mathematics. Early research efforts focused primarily on applications in control theory \cite{podlubny2002fractional} and viscoelastic materials \cite{bagley1983theoretical}, with particular emphasis on time-fractional derivatives due to their effectiveness in modeling memory and hereditary effects. More recently, the scope of fractional operators has been extended to spatial domains, where they have been applied to describe phenomena such as nonlocal responses \cite{patnaik2020ritz,di2013mechanically}, multiscale interactions \cite{ding2023fractional}, and various forms of anomalous or hybrid transport \cite{gomez2016modeling}. The fractional-order continuum mechanics framework \cite{sumelka2014fractional} governs nonlocal behavior through FDEs, modeling the mechanical response at a point as a function of field variables distributed over a finite spatial domain (nonlocal horizon). Patnaik \etal \cite{patnaik2020ritz} extended this framework to generalized definitions for the nonlocal horizon and applied to numerical studies on Euler-Bernoulli beams to explore the effect of nonlocal interactions on the elastic response. Following this, a series of studies were conducted on different structural elements \cite{sidhardh2022fractional,patnaik2021fractional}, stability of nonlocal structures \cite{sidhardh2021fractional,sidhardh2021analysis}. In these studies, the positive-definite definition following fractional-order constitutive relations for nonlocal mechanics is established, and the numerical results demonstrate that the paradoxes presented by classical nonlocal models are resolved. 
Extending the fractional calculus-based constitutive framework for nonlocal mechanics to multiphysics, Sidhardh \etal \cite{sidhardh2021thermodynamics} and Desai \etal \cite{desai2025fractional} explored long-range interactions over thermo-mechanical and piezoelectric coupling, respectively. Interestingly, enhanced electromechanical coupling is realized via tuning of the nonlocal interactions across the domain. Similarly, Malapeta \etal explored the multiscale nature of fractional-order constitutive relations for the regularization of ill-posedness and the consequent numerical inconsistencies of continuum damage mechanics~\cite{sundar2025FcdmFEM}.

Although there have been several studies on the applications of fractional calculus for multiscale and multiphysics studies, there are limited quantitative studies in the literature. The practical use of fractional-order derivatives is often limited by the mathematical complexity of their differ-integral definitions \cite{oldham1974fractional}. This complexity makes it particularly challenging to obtain analytical solutions for FDEs, especially in nonlinear cases or for general boundary and initial value problems \cite{li2011numerical}. Traditional numerical solvers designed for integer-order equations are generally unsuitable for FDEs due to weak singularities introduced by power-law kernels, which require specialized numerical treatment. To address these challenges, several numerical schemes have been developed. Finite difference methods based on the Grünwald-Letnikov definition have been applied to model viscoelastic behavior \cite{galucio2004finite, cortes2007finite}. Sweilam \etal \cite{sweilam2007numerical} implemented a semi-analytical approach based on the variational iteration method and the homotopy perturbation method to solve FDEs. Following this method, the solution is expressed as a rapidly convergent series with easily computable components in these semi-analytical schemes. Recently, Han and Lu\cite{han2024novel} applied the Fourier spectral method for numerical study on the space-variable Reisz definition for the fractional-order derivative. Their numerical experiments highlight the high efficiency and low computational cost of spectral methods for fractional-order derivatives.

More recently, finite element (FE) methods have been used for nonlocal elasticity \cite{patnaik2020ritz, sidhardh2020geometrically} and anomalous transport processes \cite{bu2014galerkin, yang2017finite}. However, FE-based solvers face limitations such as high mesh resolution requirements, globally constructed mesh distributions, or dual mesh strategies to separate local and nonlocal contributions \cite{patnaik2021fractional, patnaik2022fractional, ding2022multiscale}. Additionally, assembling system matrices becomes increasingly complex due to long-range interactions and mesh-based interpolation, particularly for non-uniform meshes.

Meshfree numerical methods, which rely on an arbitrary distribution of unconnected nodes, offer significant advantages over conventional FE approaches, particularly in handling complex and irregular domains. In contrast to the finite element method, which approximates field variables using discretized elements and local natural coordinates, meshfree methods construct the approximation functions directly in the global Cartesian coordinate system, relying solely on nodal information without the need for predefined mesh connectivity. This involves compactly supported approximation functions that form a partition of unity over the domain. In addition to this, the continuity and completeness of the approximation functions can be independently controlled. This decoupling allows for greater flexibility in the accuracy of the solution and significantly simplifies $p$-adaptivity. Moreover, specialized basis functions can be embedded to capture key features of the underlying partial differential equations (PDEs), and node distributions can be chosen in a manner to introduce discontinuities in a seamless manner. Among meshfree methods, the element-free Galerkin (EFG) method proposed by Belytschko \etal \cite{belytschko1994element} is one of the most rigorously developed and widely adopted formulations in the field. It has been successfully applied to model elastic behavior in structures governed by classical (integer order) elasticity \cite{lu1994new,krysl1995analysis,belytschko1995crack,belytschko1996dynamic} and nonlocal elasticity theories \cite{pan2010computational,zhang2016nonlocal}. More recently, Rajan \etal \cite{rajan2024element} developed a meshfree EFG approach to solve 1D FDEs. Also, Ghorbani and Semperlotti \cite{ghorbani2025one} developed the mathematical framework of the Smoothed Particle Hydrodynamics (SPH) method to numerically approximate fractional-order operators and demonstrated its efficacy over studies on nonlocal elasticity using fractional-order constitutive models.

This 1D f-EFG developed in \cite{rajan2024element} is limited to 1D fractional-order boundary value problems. Therefore, it was applied for a numerical study on nonlocal elasticity in beams. Within the 1D f-EFG method, the global nature of the Moving Least Squares (MLS) approximants is exploited for the evaluation of independent field variables and their fractional-order derivatives, where the latter are characterized by long-range, nonlocal interactions \cite{rajan2024element}. More clearly, field variables at any point are estimated using MLS approximants, utilizing the relevant nodal values from neighboring nodes within a support domain. This inherently nonlocal interpolation distinguishes the EFG approach from traditional FEM, where the approximation is limited to immediately adjacent nodes within an element. Furthermore, the f-EFG method has been shown to be computationally more efficient than equivalent mesh-based methods.

Despite these advantages, to the best of our knowledge, there has been no prior effort to explore or develop 2D-meshfree numerical solvers specifically tailored for FDEs. In contrast to the 1D domain, there are numerous additional challenges in the evaluation of 2D fractional-order derivatives. For instance, unlike the 1D fractional-order derivatives evaluated over a line, the 2D derivatives can be defined over complex and irregular domains. The suitability of meshfree methods to handle such irregular domains, compared to mesh-based approaches, is well documented\cite{belytschko1994element}. 
Moreover, challenges in numerically evaluating fractional-order derivatives limit fractional-order boundary value problems to regular geometries.    
Further, unlike 1D boundary value problems, 2D problems include mixed derivatives, and fractional-order derivatives violate Clairaut's theorem presenting $\frac{\partial^\alpha \square}{\partial x \partial y}\neq \frac{\partial^\alpha \square}{\partial y \partial x}$. 
This study aims to bridge that gap by extending the EFG method to solve 2D boundary value problems governed by fractional-order models.

The paper begins with a brief discussion of fractional-order derivatives, highlighting their relevance to elastic constitutive modeling. This theoretical background provides the basis for formulating a fractional-order boundary value problem, which serves to demonstrate the meshfree numerical method developed in this work. Although nonlocal elasticity is used as an example case to establish the governing equations, the proposed solver can be readily adapted to solve general fractional-order boundary value problems. Section~\ref{sec: review} presents the governing equations for a fractional-order nonlocal plate subject to infinitesimal and finite deformations; linear and nonlinear integro-differential equations, respectively. Section~\ref{sec: fEFG} introduces interpolation functions via the MLS approximation and outlines the numerical evaluation of their 2D fractional derivatives, a step complicated by the nonlocal nature of fractional operators. These derivatives are used to assemble the system matrices for the bending response of the plate. The efficacy of the proposed formulation is validated against benchmark solutions in Section~\ref{sec: results}.

\section{Governing equations}
\label{sec: review}
The governing equations for the nonlocal response of a fractional-order isotropic plate subject to a uniformly distributed force are derived here. Nonlocal interactions across the plate are modeled by a fractional-order constitutive model\cite{sidhardh2020geometrically}. Thereby, the elastic response of the nonlocal plate is provided by fractional-order governing differential equations. 

\subsection{Fractional-order kinematic constitutive relations}
\begin{figure}[H]
    \centering
    \begin{subfigure}{.5\textwidth}
        \centering
        \includegraphics[width=1\linewidth]{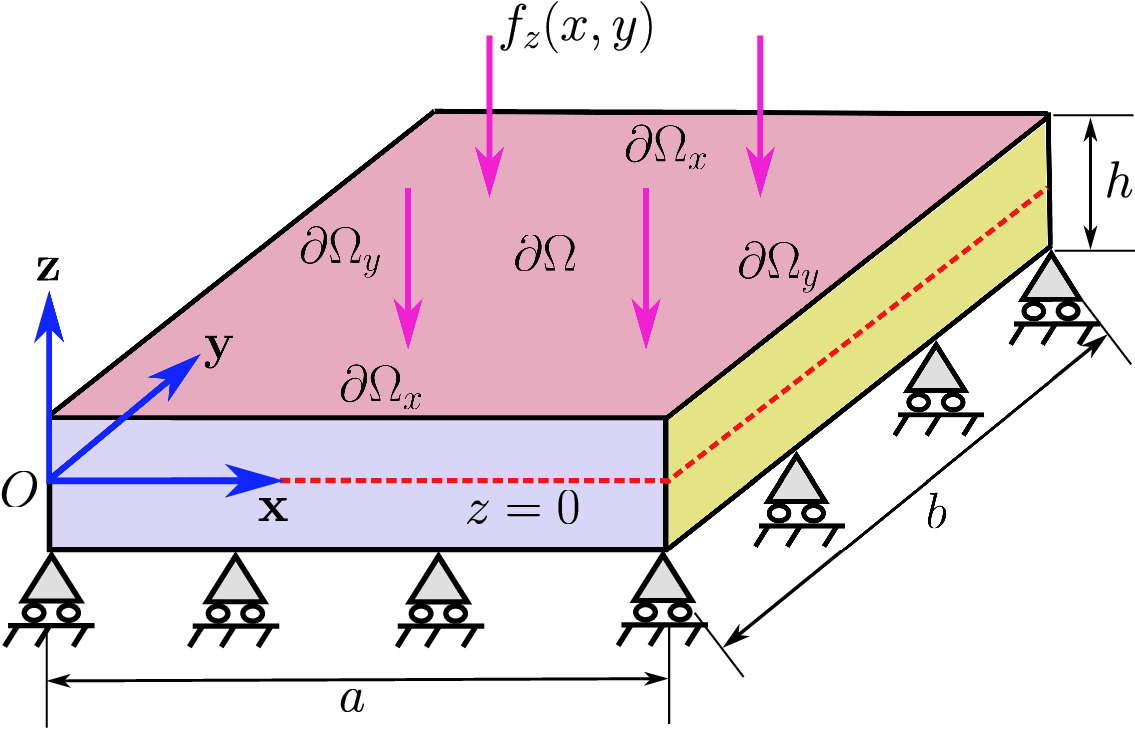}
        \caption{Rectangular plate}
        \label{fig:SSSS_plate}
    \end{subfigure}%
    \begin{subfigure}{.5\textwidth}
        \centering
        \includegraphics[width=1\linewidth]{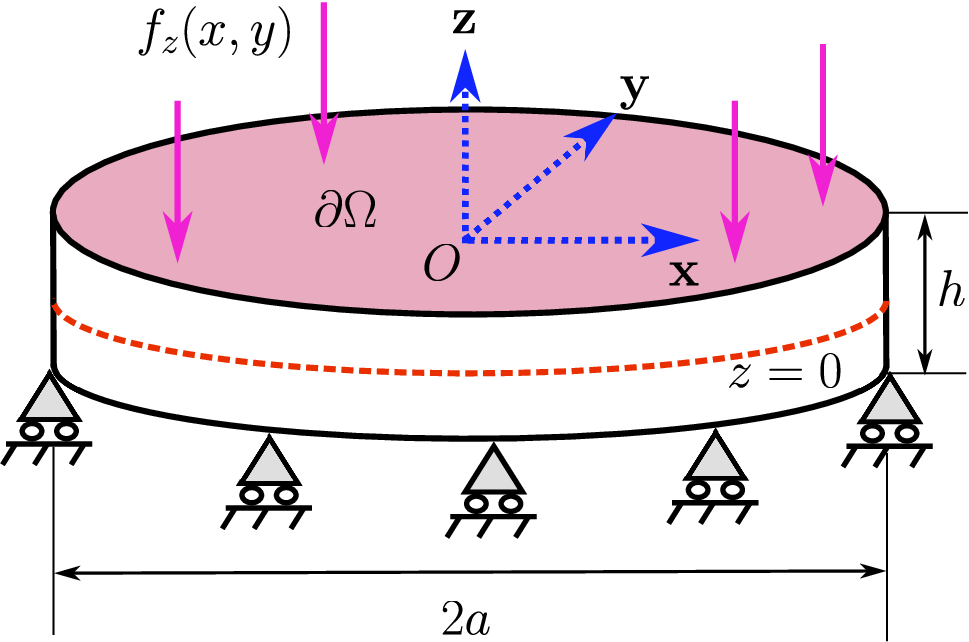}
        \caption{Circular plate}
        \label{fig:SSSS_CIRC_plate}
    \end{subfigure}
    \caption{Illustration of isotropic Kirchhoff plates with respective choices for Cartesian coordinates. Plates are subject to transverse loads on top surface.}
    \label{fig:plate_bcs}
\end{figure}

A schematic illustration of the plates chosen for the current study is illustrated in Fig.~\ref{fig:plate_bcs}. 
For the current study, two different cases of rectangular and circular plates are chosen. The in-plane dimensions of the rectangular plate are given by: $a$ and $b$, and the radius of the circular plate is considered to be $a$.
The thickness of the plate is considered along the $z$-direction, ranging from $z = -h/2$ to $z = h/2$, where $h$ is the total thickness of the plate. In this study, the plates are subjected to simply supported boundary conditions along all edges and a uniformly distributed transverse load (UDTL) acting in the $z$-direction, as illustrated in Fig.~\ref{fig:SSSS_plate}. The chosen Cartesian coordinate system is also illustrated in the figure.  For the rectangular plate, the origin is at the left-bottom corner of the mid-plane; for the circular plate, it is at the center in the mid-plane.

Assuming the plate is sufficiently thin to satisfy the classical Kirchhoff plate theory, the aspect ratio is chosen to be $a/h \geq 50$~\cite{patnaik2020geometrically}. Accordingly, the Kirchhoff kinematic assumptions are employed to study the geometrically nonlinear, nonlocal elastic behavior of the plate using a displacement-based formulation, expressed as~\cite{patnaik2020geometrically}:
\begin{subequations}
    \begin{equation}
        u(x, y, z) = u_0(x, y) - z \frac{\partial w_0(x, y)}{\partial x},
    \end{equation}
    \begin{equation}
        v(x, y, z) = v_0(x, y) - z \frac{\partial w_0(x, y)}{\partial y},
    \end{equation}
    \begin{equation}
        w(x, y, z) = w_0(x, y).
    \end{equation}
    \label{kinematics_cpt}
\end{subequations}

Here, $u_0(x,y)$ and $v_0(x,y)$ denote the in-plane displacements of the mid-plane at $z=0$ along the $x$- and $y$-directions, respectively, and $w_0(x,y)$ is the transverse displacement of the mid-plane in the $z$-direction.

In the proposed study, large deformations are assumed to present with nonlinear governing equations. However, assuming moderate strains, geometric nonlinearity can be approximated by a von Kármán nonlinear strain formulation. Furthermore, assuming that nonlocal interactions are captured via fractional-order definitions for strain-displacement relations introduced in \cite{patnaik2020ritz,sidhardh2020geometrically}, the geometrically nonlinear and nonlocal strain-displacement relations are given below:
\begin{subequations}
\label{eq:frac_strain_comp}
    \begin{equation}
    \label{eq:strain_x}
        \tilde{\epsilon}_{xx} = D^{\alpha}_{x} u_0(x,y) + \frac{1}{2} \left( D^{\alpha}_{x} w_0(x,y) \right)^2 - z~D^{\alpha}_{x} \left( \frac{\partial w_0(x,y)}{\partial x} \right),
    \end{equation}
    \begin{equation}
        \label{eq:strain_y}
        \tilde{\epsilon}_{yy} = D^{\alpha}_{y} v_0(x,y) + \frac{1}{2} \left( D^{\alpha}_{y} w_0(x,y) \right)^2 - z~D^{\alpha}_{y} \left( \frac{\partial w_0(x,y)}{\partial y} \right),
    \end{equation}
    \begin{equation}
    \label{eq:strain_xy}
        \tilde{\gamma}_{xy} = D^{\alpha}_{y} u_0(x,y) + D^{\alpha}_{x} v_0(x,y) + D^{\alpha}_{x} w_0(x,y) \, D^{\alpha}_{y} w_0(x,y)- z~\left( D^{\alpha}_{y} \left(\frac{\partial w_0(x,y)}{\partial x}\right) + D^{\alpha}_{x} \left(\frac{\partial w_0(x,y)}{\partial y}\right) \right).
    \end{equation}
\end{subequations}
The other components of strain, $\gamma_{xz} = \gamma_{yz} = \epsilon_{zz} = 0$, are assumed to be identically zero. The above strain-displacement relations are defined using fractional-order derivatives to capture long-range interactions within the plate.  

The fractional derivative $D^{\alpha}_{x_j} \psi (\mathbf{x})$ in the above expressions is defined in the Riesz-Caputo (RC) sense at the point $\mathbf{x}(x,y)$ as follows~\cite{sumelka2014thermoelasticity,patnaik2020ritz}:
\begin{equation}
    \label{eq:Frac_RC_derivative}
    D^{\alpha}_{x_j} \psi (\textbf{x}) = \frac{1}{2} \Gamma(2 - \alpha) \left[
    l_{A_j}^{\alpha-1} \, {}^C_{x_j - l_{A_j}} D^\alpha_{x_j} \psi (\textbf{x})
    - l_{B_j}^{\alpha-1} \, {}^C_{x_j} D^\alpha_{x_j + l_{B_j}} \psi (\textbf{x})
    \right],~~~j=1,2.
\end{equation}
 In this expression, $\alpha \in (0,1]$ and $(x_j - l_{A_j}, x_j + l_{B_j}) \subseteq \mathbb{R}$ are the order and terminals of the fractional derivative. Finally, the gamma function is denoted by $\Gamma (\cdot)$. The left- and right-hand Caputo derivatives defined at the point $\textbf{x}$ in the above definition are given as:
\begin{subequations}
    \label{frac_left_right caputo}
    \begin{equation}
    \label{eq:Caputo left}
    {}^C_{x_j - l_{A_j}} D^\alpha_{x_j} \psi (\textbf{x}) = \frac{1}{\Gamma(1-\alpha)} \int_{x_j - l_{A_j}}^{x_j} \frac{D^1_{x_j} \psi (\textbf{s})}{|x_j - s_j|^\alpha} \, \mathrm{d}s_j,
\end{equation}
\begin{equation}
    \label{eq:Caputo right}
    {}^C_{x_j} D^\alpha_{x_j + l_{B_j}} \psi (\textbf{x}) = \frac{1}{\Gamma(1-\alpha)} \int_{x_j}^{x_j + l_{B_j}} \frac{D^1_{x_j} \psi (\textbf{s})}{|x_j - s_j|^\alpha} \, \mathrm{d}s_j.
\end{equation}
\end{subequations}
Here $D^1_{x_j}(\cdot)$ is the first order derivative with respect to $x_j$ and $s_j$ is the coordinates of a point $\textbf{s}$ in the nonlocal horizon of the point $\textbf{x}$ along the $x_j$ direction. The degree of the fractional-order derivative $\alpha$ determines the degree of nonlocal (long-range) interactions. More clearly, a lower $\alpha$ implies a greater degree of nonlocal interactions within the same horizon defined by $l_{A_j}$ and $l_{B_j}$, while a higher $\alpha$ ($<1$) leads to more localized behavior. This is due to the slower decay in the power-law kernel weight function for the interaction between points $x_j$ and $s_j$ for a lower numerical value of fractional-order. The fractional-order derivative includes nonlocal horizons $l_{A_j}$ and $l_{B_j}$ defined in the left and right directions of the point $\textbf{x}$, respectively, along the $j-$th direction; $j=1,2$ for the 2D domain. The nonlocal horizons are defined to be a spatial variable, allowing them to be position-dependent. This helps the fractional-order framework to capture long-range interactions while being appropriately truncated at material and geometric boundaries. The integro-differential formulation of fractional-order derivatives is naturally suited to modeling long-range interactions in elastic media. Therefore, the fractional-order derivatives described above effectively capture the long-range interactions inherent in nonlocal continua. This approach is particularly advantageous for modeling the nonlocal response of Kirchhoff plates, where conventional local theories present paradoxical and inconsistent observations to describe size-dependent responses.
A detailed account of the advantages of the fractional calculus framework for nonlocal elasticity is discussed in the literature \cite{patnaik2020ritz,patnaik2022displacement}, and therefore omitted here for the sake of brevity. 
Introducing the definitions for the left and right Caputo derivatives given in Eqs.~\eqref{eq:Caputo left} and \eqref{eq:Caputo right} in the RC fractional-order derivative gives~\cite{rajan2024element,desai2025fractional}:

\begin{equation}
    \label{eq:simple_Frac_RC_derivative}
    D^{\alpha}_{x_j} \psi (\textbf{x}) = \frac{1}{2} (1 - \alpha) \left[
    l_{A_j}^{\alpha-1} \, \int_{x_j-l_{A_j}}^{x_j} \frac{D^1_{x_j} \psi (\textbf{s})}{|x_j - s_j|^\alpha} \, \mathrm{d}s_j + l_{B_j}^{\alpha-1} \, \int_{x_j}^{x_j+l_{B_j}} \frac{D^1_{x_j} \psi (\textbf{s})}{|x_j - s_j|^\alpha} \, \mathrm{d}s_j \right].
\end{equation}
The above expression can be recast in a simpler form by defining an attenuation function, $\mathcal{A}(x_j, s_j, l_{A_j},l_{B_j}, \alpha )$ to give~\cite{patnaik2020ritz}:

\begin{equation}
    \label{eq:Attenuation_simple_Frac_RC_derivative}
    D^{\alpha}_{x_j} \psi (\textbf{x}) =
     \int_{x_j-l_{A_j}}^{x_j+l_{B_j}} \, \mathcal{A}(x_j, s_j, l_{A_j}, \alpha ) ~ {D^1_{x_j} \psi (\textbf{s})} \, \mathrm{d}s_j, 
\end{equation}
where $\mathcal{A}(x_j, s_j, l_{A_j},l_{B_j}, \alpha )$ is the power-law kernel or attenuation function of the fractional-order derivative. This kernel function is analogous to the nonlocal weighting functions typically used in alternate formulations, and can be expressed as the interaction between point $x_j$ and a point $s_j$ within its horizon of nonlocal influence as~\cite{sidhardh2020geometrically}:
\begin{equation}
    \label{Attenuation_functiom}
    \mathcal{A}(x_j, s_j, l_{A_j}, l_{B_j}, \alpha) = 
    \begin{cases}
        \left( \dfrac{1 - \alpha}{2} \right) l_{A_j}^{\alpha - 1} \dfrac{1}{|x_j - s_j|^{\alpha}}, & s_j \in (x_j - l_{A_j}, x_j), \\[10pt]
        \left( \dfrac{1 - \alpha}{2} \right) l_{B_j}^{\alpha - 1} \dfrac{1}{|x_j - s_j|^{\alpha}}, & s_j \in (x_j, x_j + l_{B_j}).
    \end{cases}
\end{equation}
Note the power-law decay with respect to the relative distance $|x_j - s_j|$, modulated by the fractional order $\alpha \in (0,1]$. Clearly, as $\alpha \to 1$, the kernel approximates the local behavior, the nonlocal effects diminish, and the formulation converges to the classical elastic response of a solid body. 

It is important to note that the nonlinear strain-displacement relations incorporate fractional-order gradients of the deformation field. Specifically, the strain components $\tilde{\epsilon}_{ij}$ defined above incorporate fractional-order displacement gradients, which are assumed to be of order $\mathcal{O}(\varepsilon)$ under the moderate strain assumption:
\begin{equation}
\label{eq:linear_strain_comp}
    \left\{
        D^{\alpha}_x u_0, \;
        D^{\alpha}_y v_0, \;        
        D^{\alpha}_{x} \left( \frac{\partial w_0}{\partial x} \right), \;
        D^{\alpha}_{y} \left( \frac{\partial w_0}{\partial y} \right),
        \;
        D^{\alpha}_{y} \left( \frac{\partial w_0}{\partial x} \right), \;
        D^{\alpha}_{x} \left( \frac{\partial w_0}{\partial y} \right)
    \right\} = \mathcal{O}(\varepsilon).
\end{equation}
In addition, the formulation of geometrically nonlinear fractional-order plate theory introduces higher-order displacement gradient terms of order $\mathcal{O}(\varepsilon^2)$, such as:
\begin{equation}
\label{eq:nonlinear_strain_comp}
\left\{
    \left(D^{\alpha}_x w_0\right)^2, \;
    \left(D^{\alpha}_y w_0\right)^2, \;
    D^{\alpha}_x w_0 ~ D^{\alpha}_y w_0
\right\} = \mathcal{O}(\varepsilon^2).
\end{equation}
These nonlinear terms are essential to capture large-deformation kinematics accurately.

\subsection{Fractional-order equilibrium equations}
\label{subsec: gdes_frac}
The stress state at a material point $\mathbf{x}$ in the nonlocal plate is expressed in terms of fractional-order strain as~\cite{sidhardh2021thermodynamics}:
\begin{equation}
    \label{frac_stress}
    \tilde{\sigma}_{ij}(\mathbf{x}) = C_{ijkl} \, \tilde{\epsilon}_{kl}(\mathbf{x}), 
\end{equation}
where $C_{ijkl}$ denotes the fourth-order constitutive stiffness tensor of an isotropic solid. The resulting stress field is inherently nonlocal due to its dependence on the fractional-order strain field defined in Eq.~\eqref{eq:frac_strain_comp}. Assuming the plate is composed of isotropic material, only two independent elastic constants are employed: Young’s modulus ($E$) and Poisson’s ratio ($\nu$). However, the proposed geometrically nonlinear fractional-order plate formulation can be extended to general material symmetries using the previously defined stress and fractional-order strain fields.

The strong form of the governing differential equations for the fractional-order plate is derived using the principle of minimum potential energy: $\delta \Pi=0$\cite{reddy2006theory}. Here, $U$ is the deformation energy, and $V$ is the work done by an externally applied force. For the choice of the Cartesian coordinate system and the loading conditions illustrated in Fig.~\ref{fig:plate_bcs}, these functionals are defined as follows~\cite{reddy2015introduction}:
\begin{equation}
    \label{eq:energy_comp}
       \delta U = \int_{\partial\Omega} \int_{-h/2}^{+h/2} \tilde{\sigma}_{ij} \, \delta \tilde{\epsilon}_{ij} \, \mathrm{d}z \, \mathrm{d}S,~~~
       \delta V = \int_{\partial\Omega} \left( f_z \, \delta w_0 \right) \mathrm{d}S.
\end{equation}   
Here, $f_z$ is the UDTL on the 2D domain, $\partial \Omega$. Finally, $\mathrm{d}S = \mathrm{d}x\,\mathrm{d}y$ represents an infinitesimal area element.

By substituting Eq.~\eqref{eq:energy_comp} and \eqref{eq:frac_strain_comp} into the principle of minimum potential energy and employing the fundamental lemma of variational calculus, the corresponding governing differential equations for the in-plane response of the fractional-order plate are obtained as~\cite{patnaik2020geometrically,patnaik2021fractional}:
\begin{subequations}
    \label{eq:inplane_gde_cpt}
    \begin{equation}
       \mathcal{D}_x^{\alpha} \mathcal N_{xx} + \mathcal{D}_y^{\alpha} \mathcal N_{xy} = 0,~~~~
       \mathcal{D}_x^{\alpha} \mathcal N_{xy} + \mathcal{D}_y^{\alpha} \mathcal N_{yy} = 0.
    \end{equation}
$\{\mathcal N_{xx}, \mathcal N_{yy}, \mathcal N_{xy}\}$ are the in-plane stress resultants, defined as~\cite{reddy2015introduction}:
\begin{equation}
    \label{eq:inplane_stress_resultants}
    \{\mathcal N_{xx}, ~~\mathcal N_{yy}, ~~\mathcal N_{xy}\} = \int_{-h/2}^{+h/2} \{\tilde{\sigma}_{xx}, ~~\tilde{\sigma}_{yy}, ~~\tilde{\sigma}_{xy}\}\,\mathrm{d}z.
\end{equation}
\end{subequations}
The operator $\mathcal{D}^{\alpha}_{x_j}(\cdot)$ represents the Riesz-Riemann-Liouville (RRL) fractional derivative of order $\alpha$ with respect to the spatial coordinate $x_j$, and is defined as~\cite{patnaik2021fractional}:
\begin{equation}
    \label{eq:RRL_derivative}
        \mathcal{D}^{\alpha}_{x_j} \psi(\mathbf{x}) = \frac{1}{2} \Gamma(2 - \alpha) \left[
         l_{B_j}^{\alpha - 1} \, {}^{RL}_{x_j - l_{B_j}} D^\alpha_{x_j} \psi (\mathbf{x})
         - l_{A_j}^{\alpha - 1} \, {}^{RL}_{x_j} D^\alpha_{x_j + l_{A_j}} \psi (\mathbf{x})
    \right],~~~~j=1,2.
\end{equation}
Here, $l_{A_j}$ and $l_{B_j}$ represent the nonlocal interaction length scales defined earlier and ${}^{RL}D^\alpha$ denotes the Riemann-Liouville fractional derivative of order $\alpha$, defined over the interval $(x_j - l_{B_j},\, x_j + l_{A_j})$. The function $\psi(\mathbf{x})$ is an arbitrary field, and the expressions ${}^{RL}_{x_j - l_{B_j}} D^{\alpha}_{x_j} \, \psi(x)$ and ${}^{RL}_{x_j} D^{\alpha}_{x_j + l_{A_j}} \, \psi(x)$ correspond to the left- and right-sided Riemann-Liouville fractional derivatives of $\psi(x)$, respectively. Note the integro-differential nature of the above 2D governing equations. More clearly, the fractional-order derivatives (both RC and RRL) in the governing equations given above render them to be integro-differential. 

The corresponding boundary conditions for the fractional-order Kirchhoff plate on the edges $\partial \Omega_x$ and $\partial \Omega_y$ illustrated in Fig.~\ref{fig:plate_bcs} are prescribed as follows:~\cite{patnaik2020geometrically,patnaik2021fractional}:
\begin{equation}
    \label{eq:Bcs_inplane}
    \forall \textbf{x} ~ \text{on}~\partial \Omega_x \cup \partial\Omega_y: \quad \delta u_0=\delta v_0=0.
\end{equation}
Recall that a simply supported boundary condition is considered in the present study. Accordingly, the essential boundary conditions for the in-plane variables are imposed.

The governing differential equation for the nonlinear transverse response of the fractional-order Kirchhoff plate is expressed as~\cite{patnaik2020geometrically,patnaik2021fractional}:
\begin{subequations}
    \label{eq:transverse_GDE}
 \begin{equation}
    \begin{aligned}
        & {D}_x^{1} \left[ \mathcal{D}_x^{\alpha} \mathcal M_{xx} \right]
        + {D}_x^{1} \left[ \mathcal{D}_y^{\alpha} \mathcal M_{xy} \right]
        + {D}_y^{1} \left[ \mathcal{D}_y^{\alpha} \mathcal M_{yy} \right]
        + {D}_y^{1} \left[ \mathcal{D}_x^{\alpha} \mathcal M_{xy} \right] \\
        & + \mathcal{D}_x^{\alpha} \left( \mathcal N_{xx} {D}_x^{\alpha} w_0 + \mathcal N_{xy} {D}_y^{\alpha} w_0 \right)
        + \mathcal{D}_y^{\alpha} \left( \mathcal N_{xy} {D}_x^{\alpha} w_0 + \mathcal N_{yy} {D}_y^{\alpha} w_0 \right)
        + f_z = 0.
\end{aligned}
\end{equation}
$\{\mathcal M_{xx},~~\mathcal M_{yy},~~\mathcal M_{xy}\}$ are the moment resultants, defined as~\cite{reddy2015introduction}:
\begin{equation}
    \label{eq:outplane_stress_resultants}
        \{\mathcal M_{xx}, \mathcal M_{yy}, \mathcal M_{xy}\} = \int_{-h/2}^{+h/2} 
\left\{ z\tilde{\sigma}_{xx},\; z\tilde{\sigma}_{yy},\; z\tilde{\sigma}_{xy} \right\} 
\, \mathrm{d}z.
\end{equation}   
\end{subequations}
Recall that the stresses in the above definitions are expressed in terms of fractional-order strains following material constitutive relations given in Eq.~\eqref{frac_stress}.

The essential and natural boundary conditions corresponding to the transverse displacement are expressed as follows~\cite{patnaik2020geometrically,patnaik2021fractional}:
\begin{subequations}\label{eq:Bcs}
    \begin{equation}
        \forall \textbf{x} ~ \text{on}~\partial \Omega_x \cup \partial\Omega_y: \quad \delta w_0=0,
    \end{equation}
    and
    \begin{align}
       \forall \textbf{x} ~ \text{on}~ \partial \Omega_x:~
            \delta D_{x}^1 w_0=0\quad\text{or}\quad I_x^{1 - \alpha}\mathcal M_{xx} = 0,\\
        \forall \textbf{x} ~ \text{on} ~ \partial \Omega_y:~
            \delta D_{y}^1 w_0=0\quad\text{or}\quad I_y^{1 - \alpha}\mathcal M_{yy} = 0.
    \end{align}
\end{subequations}
where, the Riesz fractional integral $I_{x_j}^{1 - \alpha}(\cdot)$ is defined as~\cite{patnaik2020geometrically,patnaik2021fractional}:
\begin{equation}
\label{eq:Reisz_integral}
 I_{x_j}^{1-\alpha} \psi(\textbf{x}) = \frac{1}{2} \Gamma(2 - \alpha) \left[
    l_{B_j}^{\alpha - 1} \, {}_{x_j - l_{B_j}} I^{1-\alpha}_{x_j} \psi (\mathbf{x})
    -  l_{A_j}^{\alpha - 1} \, {}_{x_j} I^{1-\alpha}_{x_j+l_{A_j}} \psi (\mathbf{x})
    \right],~~~j=1,2.
\end{equation}
The Riesz fractional integral $I^{1-\alpha}_{x_j} \psi(\mathbf{x})$ is defined in terms of fractional-order Reimann integrals over the interval $(x_j - l_{B_j},\, x_j + l_{A_j})$. Specifically, the term ${}_{x_j - l_{B_j}} I^{1-\alpha}_{x_j} \psi(\mathbf{x})$ corresponds to the left-sided fractional integral, while ${}_{x_j} I^{1-\alpha}_{x_j + l_{A_j}} \psi(\mathbf{x})$ represents its right-sided counterpart, both of order $(1 - \alpha)$. 

The essential boundary condition corresponding to the transverse displacement and the natural boundary conditions corresponding to the moments are applicable for the simply supported boundary condition considered in this study.
The 2D FDEs given in Eqs.~\eqref{eq:inplane_gde_cpt} and \eqref{eq:transverse_GDE} determine the geometrically nonlinear response of nonlocal Kirchhoff rectangular plates. These fractional-order PDEs, along with the corresponding BCs given in Eqs.~\eqref{eq:Bcs_inplane} and \eqref{eq:Bcs} require to be numerically solved for a quantitative analysis. Although the above equations are developed for the rectangular plate, they can be similarly extended for the circular plate. The derivation of the same is omitted here for the sake of brevity.

In the following section, we present a detailed numerical framework developed using the 2D-EFG (meshfree) method for the solution of these FDEs for geometrically nonlinear response. It is important to note that the nonlinear FDEs for the nonlocal Kirchhoff plate can be reduced to a linear form by considering only the linear ($\mathcal{O}(\varepsilon)$) fractional-order displacement gradient terms in Eq.~\eqref{eq:linear_strain_comp} of the strain-displacement relations. More clearly, the nonlinear terms of order $\mathcal{O}(\epsilon^2)$ in the strain-displacement relations given in Eq.~\eqref{eq:nonlinear_strain_comp} are neglected for the linear analysis of plates. The subsequent material constitutive relations and the governing equations are updated accordingly. This linearization is relevant for an infinitesimal deformation of the nonlocal plates.

\section{Element-free Galerkin model for fractional-order plate}
\label{sec: fEFG}
A numerical framework using EFG method is developed in this section for the previously derived 2D FDEs.
The proposed f-EFG formulation is developed to efficiently capture the nonlocal behavior of fractional-order Kirchhoff plate models.
The meshfree nature of the formulation makes it particularly advantageous for fractional-order problems that involve differ-integral operators, as it facilitates the modeling of long-range interactions through a flexible distribution of nodes rather than relying on fixed element connectivity. Owing to this versatility, the developed solver can be readily applied to a wide range of 2D FDEs involving spatial fractional derivatives.

The framework comprises three key stages: (i) formulation of MLS interpolation functions for constructing smooth and differentiable approximations over scattered nodes, (ii) extension of these MLS approximants to compute fractional-order derivatives based on fractional-order derivative operators, and (iii) derivation of the weak form and algebraic system matrices used in both linear and nonlinear response of the Kirchhoff plates. 

\subsection{MLS approximation functions}
\label{subsec: mls_func}
The construction of MLS approximation functions, in the context of integer-order differential equations, is well established~\cite{belytschko1994element}. Furthermore, it has previously been applied to fractional-order 1D PDEs~\cite{rajan2024element}. 
However, it is worth noting that applying 1D MLS approximants to fractional-order PDEs is relatively straightforward due to the inherent alignment between the nonlocal horizon of influence and the support domain of the weight function. In contrast, the present study extends this formulation to two dimensions, utilizing 2D basis functions along with a weight function defined over a 2D support domain. The introduction of 2D fractional-order derivatives leads to a nonlocal horizon of influence in two spatial directions, which may differ from (or coincide with) the chosen support domain.
In the current context, the MLS approximants present smooth and high-order interpolations of unknown field variables over a set of arbitrarily distributed nodes in the problem domain. This is unlike conventional finite element methods that rely on predefined connectivity and $C^1$-continuous shape functions for an accurate interpolation of the higher-order derivatives in the formulation for Kirchhoff plates. Further to this, the flexibility in node distribution for interpolation for EFG is well-suited for application in problems for nonlocal interactions. The first advantage applies to the current example of a Kirchhoff plate, whereas the second is particularly suited for the nonlocal formulation modeled here using fractional-order derivatives.

A brief overview of the steps involved in the derivation of the MLS approximants for the unknown field variable $\mathrm{w}_{0}(\mathbf{x})$, where \( \mathbf{x} = \{x, y\}^T \) (with \( T \) denoting the transpose), is presented here. For a comprehensive mathematical treatment, the reader is referred to the work of Lancaster and Salkauskas in~\cite{lancaster1981surfaces}. Consider the 2D spatial domain, corresponding to the mid-plane of the plate, defined as \( x \in [0, a] \) and \( y \in [0, b] \) for the rectangular plate and \( x \in [-a, a] \) and \( y \in [-a, a] \) for the circular plate, discretized using \( n \)-nodes. For convenience, unless otherwise mentioned, this study considers a uniform distribution of nodes. At any arbitrary point \( \mathbf{x} \) within the domain, the local approximation of the unknown field variable is constructed using the MLS interpolation as~\cite{belytschko1994element,liu2005meshfree,atluri1999mlpg}:
\begin{subequations}
\begin{equation}
\mathrm{w}_{0}(\mathbf{x}) \approx  \left\{ \mathbf{p}(\mathbf{x}) \right\}^T \, \left\{\mathbf{a}(\mathbf{x})\right\}.
\label{eq:mls_approx}
\end{equation}
In the above expression,
\begin{equation}
\left\{ \mathbf{p}(\mathbf{x}) \right\}^T = [1~~x~~y~~x^2~~xy~~y^2] ~~~~(m=6),
\label{eq:quadratic_basis}
\end{equation}
is the polynomial basis vector, and \( m \) denotes the number of terms in the polynomial basis vector. In the present study, a quadratic monomial basis is adopted for the sake of completeness and continuity of derived approximations\cite{belytschko1994element}. Further, 
\begin{equation}
\left\{\mathbf{a}(\mathbf{x})\right\}^T = [a_0(\mathbf{x}) ~~a_1(\mathbf{x})~~\cdots~~a_5(\mathbf{x})],
\label{eq:coeff_vector}
\end{equation}
is the vector of unknown coefficients. These coefficients are evaluated by minimizing a weighted least-squares functional defined over a neighborhood around \( \mathbf{x} \), thereby ensuring a smooth and mesh-independent approximation of the field.
\end{subequations}

More clearly, a weighted \( \mathbb{L}_2 \)-norm of the approximation error is defined and minimized over the \( n \) nodes. This \( \mathbb{L}_2 \)-norm is defined as follows:
\begin{equation}
J(\left\{\mathbf{a}(\mathbf{x})\right\}) = \sum_{i=1}^{n} \widehat{W}(\mathbf{x} - \mathbf{x}^i) 
\left[ \left\{\mathbf{p}(\mathbf{x}^i)\right \}^T \, \left\{\mathbf{a}(\mathbf{x})\right\} - w_0(\mathbf{x}^i) \right]^2.
\label{eq:least_squares}
\end{equation}
In the above expression, \( \mathbf{x}^i \) denotes the coordinates of the \( i^{\text{th}} \) node, and \(\left\{ \mathbf{p}(\mathbf{x}^i)\right\} \) represents the polynomial basis vector evaluated at that node. \( \widehat{W}(\mathbf{x} - \mathbf{x}^i) \) is a compactly supported positive-definite weight function centered at the evaluation point \( \mathbf{x} \). The choice of this weight function and the corresponding support domain plays a crucial role in enforcing the desired smoothness, locality, and continuity characteristics required for the MLS approximation functions.

In the present study, the cubic spline weight function is chosen due to the $C^1$ continuity requirement of the approximations~\cite{liu2005meshfree}. The weight function defined along $x$-direction is given by:
\begin{equation}
\widehat{W}(x - x^i) = \widehat{W}(r_x^i) =
\begin{cases}
\frac{2}{3} - 4r^2 + 4r^3, & \text{if } r \leq 0.5, \\
\frac{4}{3} - 4r + 4r^2 - \frac{4}{3}r^3, & \text{if } 0.5 < r \leq 1, \\
0, & \text{if } r > 1,
\end{cases}
\label{eq:cubic_weight}
\end{equation}
where \( r_x^i \) is the normalized distance along $x-$direction, between the point \( \mathbf{x}(x,y) \) and the \( i^{\text{th}} \) node \(\mathbf{x}^i( x^i ,  y^i)\) in its support domain. More specifically, normalized distance variables along the two coordinate axes for point $\mathbf{x}$ w.r.to point $\mathbf{x}^i$ are defined as:
\begin{equation}
r_x^i = \frac{|x - x^i|}{d_{m,x}}~~~\text{and}~~~\quad r_y^i = \frac{|y - y^i|}{d_{m,y}}.
\label{eq:normalized_distance}
\end{equation}
Here, \( d_{m,x} \) and \( d_{m,y} \) represent the length of the support domain in the \( x \)- and \( y \)-directions, respectively. 
The subscript refers to the spatial direction of evaluation (e.g. \( x \) or \( y \)), while the superscript \( i \) indicates the evaluation at the \( i^{\text{th}} \) node. The weight function for the 2D domain is defined as a tensor product of 1D weight functions in each coordinate direction, yielding~\cite{belytschko1994element,dolbow1998introduction}:
\begin{equation}
\widehat{W}(\mathbf{x} - \mathbf{x}^i) = \widehat{W}(r_x^i) \cdot \widehat{W}(r_y^i).
\label{eq:tensor_weight}
\end{equation}
The error norm in Eq.~\eqref{eq:least_squares} is recast to be:
\begin{subequations}
\begin{equation}
J(\left\{a(\mathbf x)\right\}) = \big[ [P]\left\{\mathbf{a}(\mathbf{x})\right\} - \{\bm {w}_0\} \big]^T~\big[W\big]~\big[ [P]\left\{\mathbf{a}(\mathbf{x})\right\} - \{\bm {w}_0\} \big],
\label{eq:L2_error}
\end{equation}
where, the matrices $[P]$ and $[W]$ consolidate the polynomial basis functions and weight function evaluated at $n$-nodes following their definitions in Eq.~\eqref{eq:quadratic_basis} and Eq.~\eqref{eq:tensor_weight}, respectively. The expressions for each of these matrices are given below:
\begin{equation}
[P] =
\begin{bmatrix}
\left \{\mathbf{p}^T(\mathbf{x}^1)\right \} \\
\left \{\mathbf{p}^T(\mathbf{x}^2)\right \} \\
\vdots  \\
\left \{\mathbf{p}^T(\mathbf{x}^n)\right \}
\end{bmatrix}_{n \times m}
~~\text{and}~~
[W] =
\begin{bmatrix}
\widehat{W}(x - x^1) & 0 & \cdots & 0 \\
0 & \widehat{W}(x - x^2) & \cdots & 0 \\
\vdots & \vdots & \ddots & \vdots \\
0 & 0 & \cdots & \widehat{W}(x - x^n)
\end{bmatrix}_{n \times n}.
\label{eq:P_W_matrices}
\end{equation}
In the above equation, $\{\bm {w}_0\}$ is the vector of unknown nodal values of the transverse displacement field variable, defined as:
\begin{equation}
\{\bm {w}_0\}^T = 
\begin{bmatrix}
w_0(x^1) \quad
w_0(x^2) \quad
\cdots \quad
w_0(x^n)
\end{bmatrix}_{1 \times n}.
\label{eq:w0_vector_transpose}
\end{equation}
\end{subequations}
Finally, the column vector $\left\{\mathbf{a}(\mathbf{x})\right\}$ in the above equation is defined in Eq.~\eqref{eq:coeff_vector}. To determine this unknown vector, the error functional \( J \) is minimized by enforcing the stationary condition with respect to $\left\{\mathbf{a}(\mathbf{x})\right\}$:
\begin{equation}
\frac{\partial J(\{a(\mathbf x)\})}{\partial \left\{\mathbf{a}(\mathbf{x})\right\}} = [A(\mathbf x)] \left\{\mathbf{a}(\mathbf{x})\right\} - [H(\mathbf x)] \{\bm {w}_0\} = 0.
\label{eq:stationarity}
\end{equation}
Solving the above equation yields:
\begin{equation}
\left\{\mathbf{a}(\mathbf{x})\right\} = [A(\mathbf x)]^{-1} [H(\mathbf x)] \{\bm {w}_0\}.
\label{eq:a_solution}
\end{equation}
In the above, the moment matrix \( [A(\mathbf x)] \) and the projection matrix \( [H(\mathbf x)] \) are given by:
\begin{equation}
[A(\mathbf x)] = 
{[P]^T} \,
[W] \,
{[P]}, 
\quad 
[H(\mathbf x)] = 
{[P]^T} \,
{[W]}.
\label{eq:A_H_definitions}
\end{equation}
Substituting the matrix definitions into the result for $\left\{\mathbf{a}(\mathbf{x})\right\}$, the approximation in Eq.~\eqref{eq:mls_approx} can be rewritten as follows~\cite{belytschko1994element,dolbow1998introduction}:
\begin{equation}
\mathrm{w}_{0}(\mathbf{x}) = 
{\left\{ \mathbf{p}(\mathbf{x}) \right\}^T}\, 
{\left\{\mathbf{a}(\mathbf{x})\right\}} =
\underbrace{
{\left\{ \mathbf{p}(\mathbf{x}) \right\}^T} \, 
{[{A}(\mathbf{x})]^{-1}} \, 
{[{H}(\mathbf{x})]}
}_{\text{approximation function } [\phi(\mathbf{x})]} 
{\{\bm {w}_0\}}.
\label{eq:MLS_approximation}
\end{equation}
The above formulation expresses the transverse displacement at any point $\mathbf{x}$ as a function of its associated nodal values. This relationship, established within a meshfree framework, is known as the MLS approximation. Specifically, the function $[\phi(\mathbf{x})]$ represents the MLS shape function (or approximant), which defines the deformation field at a point $\mathbf{x}$ in the 2D domain based on the contributions of its surrounding nodal coordinates. Recall that the FE numerical solver for Kirchhoff plate theory requires the inclusion of additional degrees of freedom to satisfy the $C^1$ continuity requirement, leading to computational complexity. In contrast, the current approach restricts the nodal coordinates to transverse displacement, and continuity is achieved by incorporating higher-order terms into the polynomial basis function and an appropriate choice for the weight function. Finally, the same procedure is used to derive the MLS approximations for in-plane field variables. For this purpose, the polynomial basis function in Eq.~\eqref{eq:quadratic_basis} is restricted to a linear basis, satisfying the \(C^0\) continuity requirement while using the same weight function, thereby producing smooth results for the in-plane field variables. 

Unlike standard finite element interpolations, the MLS approximation functions do not satisfy the Kronecker delta property. Therefore, in general, $\mathrm{w}_{0}(\mathbf x^i) \neq w_0(\mathbf x^i)$. Consequently, the approximation functions do not pass through the nodal values exactly, making the direct imposition of essential (Dirichlet) boundary conditions more difficult. To address this, the collocation technique is used to enforce the essential boundary constraints.

Furthermore, the MLS interpolant is global in nature. The approximate value \( \mathrm{w}_{0}(\mathbf x) \) at a point \( \mathbf{x} \) depends on all nodal values within the support domain of that point, as opposed to classical FEM where the approximation at a point depends only on the nodal values of the element containing that point.
This dependence on nodal values over the support domain eliminates the need for element-level connectivity, which is particularly advantageous when dealing with nonlocal formulations.  
More clearly, a key drawback of FEM in the context of fractional-order derivative evaluation at a point is the need for identification and bookkeeping of all neighboring elements that fall within the nonlocal support domain. This disadvantage becomes particularly severe for 2D fractional-order derivatives, where independent nonlocal horizon along $x-$ and $y-$axes requires separate bookkeeping. Moreover, the complexity is further increased for material points near boundaries because of a truncation of the nonlocal horizon of influence. The computation involves accumulating numerical contributions from these elements, which in turn demands additional data structures and element connectivity tracking. This increases computational complexity and leads to significant overhead, especially in high-dimensional or adaptively refined meshes. In contrast, as will be shown in the section below, the MLS approximants allow a direct application of the integro-differential operators for the numerical evaluation of fractional-order derivatives. This avoids the need for explicit element connectivity, bookkeeping, or repeated element mapping at each evaluation point, thereby reducing both algorithmic complexity and overall computational cost. A detailed account of the advantages of MLS approximants for fractional-order PDEs, albeit for 1D problems, has been discussed in \cite{rajan2024element}. The same can be extended to the current study on 2D numerical solver for fractional-order PDEs, with independent horizon of nonlocal influence along the two directions.

\subsection{MLS approximation of fractional-order derivative}
As discussed in Section~\ref{sec: review}, the nonlocal behavior of a Kirchhoff plate is captured by introducing fractional-order derivatives into the constitutive relations. For this, the strain-displacement equations are reformulated in terms of these fractional derivatives, leading to f-PDEs. Therefore, it is necessary to develop numerical approximations for an accurate evaluation of these derivatives using the MLS approximation functions formulated earlier.

The local (integer-order) derivatives of the displacement field \( w_0(\mathbf{x}) \), approximated by the MLS approximation functions \( \phi_i(\mathbf{x}) \), are expressed as follows~\cite{rajan2024element,belytschko2000efg}:
\begin{subequations}
    \label{local_derivative_approx}
    \begin{equation}
        D_x^1 w_0(\mathbf{x}) = [B_{w_0,x}(\mathbf{x})]\{\bm {w}_0\},~~~~D_y^1 w_0(\mathbf{x}) = [B_{w_0,y}(\mathbf{x})]\{\bm {w}_0\},
        \end{equation}
        \begin{equation}
        D_x^2 w_0(\mathbf{x}) = [B_{w_0,xx}(\mathbf{x})]\{\bm {w}_0\},~~~~D_y^2 w_0(\mathbf{x}) = [B_{w_0,yy}(\mathbf{x})]\{\bm {w}_0\},~~~~D_x^1 D_y^1 w_0(\mathbf{x}) = [B_{w_0,xy}(\mathbf{x})] \{\bm {w}_0\}.
    \end{equation}
   The strain-displacement matrices $[B_\square]$ corresponding to integer-order derivatives of the displacement field variable $\{\bm {w}_0\}$ can be expressed in terms of the derivatives of the MLS approximation functions derived above as follows:
    \begin{equation}
        [B_{w_0,x}(\mathbf{x})] = \dfrac{\partial }{\partial x}\left\{\phi(\mathbf{x})\right\},~~~~[B_{w_0,y}(\mathbf{x})] = \dfrac{\partial }{\partial y}\left\{\phi(\mathbf{x})\right\},
    \end{equation}
    \begin{equation}
        [B_{w_0,xx}(\mathbf{x})] = \dfrac{\partial^2 }{\partial x^2}\left\{\phi(\mathbf{x})\right\},~~~~[B_{w_0,xy}(\mathbf{x})] = \dfrac{\partial^2 }{\partial x\partial y}\left\{\phi(\mathbf{x})\right\},~~~~[B_{w_0,yy}(\mathbf{x})] = \dfrac{\partial^2 }{\partial y^2}\left\{\phi(\mathbf{x})\right\}.
    \end{equation}
\end{subequations}
The derivatives of the MLS approximation function \( \left \{ \phi(\mathbf x) \right \} \) obtained from the definition given in Eq.~\eqref{eq:MLS_approximation} are provided in Appendix. A similar definition of nonlocal $[B_\square]$ is essential to numerically approximate the fractional-order derivatives of the displacement variable. 

From the results provided in the Appendix, it is clear that the integer-order derivatives of the MLS approximants are evaluated following the classical Leibniz rule (also referred to as product rule for derivatives). However, in contrast to derivatives of integer order, derivatives of fractional-order do not satisfy the classical Leibniz rule\cite{podlubny1994fractional,tarasov2013no}. Therefore, closed-form expressions for fractional derivatives of approximation functions cannot be directly derived~\cite{tarasov2013no}. 
We implement numerical integration methods for singular functions for an approximation of the differ-integral definition of a fractional-order derivative. The Riesz-Caputo fractional-order derivative of the transverse displacement coordinate \( w_0(\mathbf x) \) is expressed using the definition in Eqs.~\eqref{eq:Frac_RC_derivative} and~\eqref{frac_left_right caputo}  as~\cite{rajan2024element}:
\begin{equation}
\begin{aligned}
D_{x_j}^\alpha w_0(\mathbf x) = \frac{1}{2}(1 - \alpha) \left( {l}_{A_j}^{\alpha-1} \int_{x_j - l_{A_j}}^{x_j} \frac{D_{x_j}^1 w_0(\mathbf s)}{|x_j - s_j|^\alpha} \, \mathrm{d}s_j \right)
+ \frac{1}{2}(1 - \alpha) \left({l}_{B_j}^{\alpha-1} \int_{x_j}^{x_j + l_{B_j}} \frac{D_{x_j}^1 w_0(\mathbf s)}{|s_j - x_j|^\alpha} \, \mathrm{d}s_j \right),~~j=1,2.
\end{aligned}
\label{eq:riesz_caputo_frac}
\end{equation}
It is important to note that the above expression contains the integer-order derivative of the displacement field. By substituting the first-order derivative approximation based on the MLS formulation given in Eq.~\eqref{local_derivative_approx}, the expression can be rewritten as:
\begin{equation}
\begin{aligned}
D_{x_j}^\alpha w_0(\mathbf x) &= \frac{1}{2} (1 - \alpha) \left[\, l_{A_j}^{\alpha - 1} \int_{x_j - l_{A_j}}^{x_j} \frac{[B_{w_0,x_j}(\mathbf s)]}{|x_j - s_j|^\alpha} \, \mathrm{d}s_j +  \, l_{B_j}^{\alpha - 1} \int_{x_j}^{x_j + l_{B_j}} \frac{[B_{w_0,x_j}(\mathbf s)]}{|x_j - s_j|^\alpha} \, \mathrm{d}s_j\right]\{\bm {w}_0\}.
\end{aligned}
\label{eq:frac_deriv_Bform}
\end{equation}
This expression can be recast as
\begin{subequations}
\label{eq:frac_deriv_compact}
    \begin{equation}
D_{x_j}^\alpha w_0(\mathbf x) = \left[\tilde{B}_{w_0,x_j}(x_j, l_{A_j},l_{B_j}, \alpha)\right]\{\bm {w}_0\},
\end{equation}
where
\begin{equation}
\left[\tilde{B}_{w_0,x_j}(x_j, l_{A_j},l_{B_j}, \alpha)\right]=\left[
\tilde{B}_{w_0,x_j}^L(x_j, l_{A_j}, \alpha)
\right]  +
\left[
\tilde{B}_{w_0,x_j}^R(x_j, l_{B_j}, \alpha)
\right].
\end{equation} 
\end{subequations}
The above definition resembles the approximation of integer-order derivatives in Eq.~\eqref{local_derivative_approx}. In the above result, $\left[\tilde{B}_{w_0,x_j}(x_j, l_{A_j},l_{B_j}, \alpha)\right]$ is the analogously defined strain-displacement matrix corresponding to the fractional-order derivative of the transverse displacement. The above matrix derives contributions from the left and right horizons of nonlocal influence, captured by:
\begin{subequations}
    \label{left_right B tildes}
    \begin{equation}
[\tilde{B}_{w_0,x_j}^L(x_j, l_{A_j}, \alpha)] = \frac{(1 - \alpha)\, 2^{\alpha - 1}}{l_{A_j}} \int_{-1}^{+1} \frac{[B_{w_0,x_j}(\zeta)]}{(1 - \zeta)^\alpha} \, \mathrm{d}\zeta,
\label{eq:B1_left}
\end{equation}
\begin{equation}
[\tilde{B}_{w_0,x_j}^R(x_j, l_{B_j}, \alpha)] = \frac{(1 - \alpha)\, 2^{\alpha - 1}}{l_{B_j}} \int_{-1}^{+1} \frac{[B_{w_0,x_j}(\zeta)]}{(1 + \zeta)^\alpha} \, \mathrm{d}\zeta.
\label{eq:B1_right}
\end{equation}
\end{subequations}
These expressions are derived from Eq.~\eqref{eq:frac_deriv_Bform} by applying a change in the integration variable, which transforms the integration limits to a form that is better suited for numerical integration. 
Recall that the MLS approximants are global in nature, which renders the left- and right-order fractional-order derivative matrices, \([\tilde{B}_{w_0,x_j}^L(x_j, l_{A_j}, \alpha)]\) and \([\tilde{B}_{w_0,x_j}^R(x_j, l_{B_j}, \alpha)]\), globally to be defined. The definitions of these matrices include integration limits of \( (x_j - l_{A_j},\, x_j) \) and \( (x_j,\, x_j + l_{B_j}) \) corresponding to respective nonlocal horizons of influence. Thus, nonlocal contributions to the left and right of \( x_j \) are trivially evaluated, and the f-EFG approach eliminates the need for additional bookkeeping to track the nonlocal influence horizon in the numerical estimation of fractional-order derivatives. 
We note here that weakly singular integrals that arise in the transformed interval \( \zeta \in [-1, 1] \) can be accurately captured using the Gauss-Jacobi quadrature rule~\cite{hale2013fast}, provided that a sufficient number of Gauss-Jacobi quadrature points are used to ensure accurate resolution. Thus, the fractional $[\tilde{B}]$ matrix connecting the transverse deformation and its gradients can be numerically evaluated.

Similarly, the other fractional-order derivatives of the mid-plane displacement coordinates can be evaluated. Some of these are provided below:
\begin{subequations}
\label{eq:frac_derivatives}
\begin{equation}
D_x^\alpha u_0(\mathbf{x}) = [\tilde{B}_{u_0,x}(\mathbf{x}, l_{A_1}, l_{B_1}, \alpha)]\, \{\bm {u}_0\},~~~D_y^\alpha u_0(x) = [\tilde{B}_{u_0,y}(\mathbf{x}, l_{A_2}, l_{B_2}, \alpha)]\, \{\bm {u}_0\},
\label{eq:frac_du}
\end{equation}
\begin{equation}
D_x^\alpha v_0(\mathbf{x}) = [\tilde{B}_{v_0,x}(\mathbf{x}, l_{A_1}, l_{B_1}, \alpha)]\, \{\bm {v}_0\},~~~D_y^\alpha v_0(x) = [\tilde{B}_{v_0,y}(\mathbf{x}, l_{A_2}, l_{B_2}, \alpha)]\, \{\bm {v}_0\},
\label{eq:frac_dv}
\end{equation}
\begin{equation}
D_x^\alpha \left( \frac{\partial w_0(\mathbf{x})}{\partial x} \right) = 
[\tilde{B}_{w_0,xx}(\mathbf{x}, l_{A_1}, l_{B_1}, \alpha)]\, \{\bm {w}_0\},~~~D_x^\alpha \left( \frac{\partial w_0(\mathbf{x})}{\partial y} \right) = 
[\tilde{B}_{w_0,yx}(\mathbf{x}, l_{A_1}, l_{B_1}, \alpha)]\, \{\bm {w}_0\},
\label{eq:frac_dxalpha}
\end{equation}
\begin{equation}
D_y^\alpha \left( \frac{\partial w_0(\mathbf{x})}{\partial x} \right) = 
[\tilde{B}_{w_0,xy}(\mathbf{x}, l_{A_2}, l_{B_2}, \alpha)]\, \{\bm {w}_0\},~~~D_y^\alpha \left( \frac{\partial w_0(\mathbf{x})}{\partial y} \right) = 
[\tilde{B}_{w_0,yy}(\mathbf{x}, l_{A_2}, l_{B_2}, \alpha)]\, \{\bm {w}_0\}.
\label{eq:frac_dyalpha}
\end{equation}
\end{subequations}
Here, \( \{\bm {u}_0\} \), \( \{\bm {v}_0\} \), and \( \{\bm {w}_0\} \) are the global vectors of the nodal displacement vectors for the in-plane and transverse displacements of the mid-plane, respectively. Note that separate approximation functions are employed for the in-plane and transverse displacement components. This distinction arises from the careful selection of the basis functions to satisfy different continuity requirements for each variable. Furthermore, unlike integer-order derivatives, where Clairaut's theorem is applied to ensure $[B_{w_0,xy}]=[B_{w_0,yx}]$, the mixed derivatives for the nonlocal plate are sensitive to the order of fractional-order and integer-order derivatives. More clearly, the $[\tilde{B}_{w_0,xy}(\mathbf{x}, l_{A_2}, l_{B_2}, \alpha)]\neq [\tilde{B}_{w_0,yx}(\mathbf{x}, l_{A_1}, l_{B_1}, \alpha)]$. This is because in the former the fractional-order derivative is along the $y-$ direction and in the latter term it is along the $x-$ direction. 

The fractional $[\tilde{B}]$ matrices corresponding to the higher-order derivatives of the transverse displacement are given as follows: 
\begin{subequations}
    \label{b2 b3 tildes}
    \begin{equation}
        [\tilde{B}_{w_0,xx}(\mathbf{x}, l_{A_1}, l_{B_1}, \alpha)] = 
        \frac{(1 - \alpha)\, 2^{\alpha - 1}}{l_{A_1}} 
        \int_{-1}^{+1} \frac{[B_{w_0,xx}(\zeta)]}{(1 - \zeta)^{\alpha}} \, \mathrm{d}\zeta +
        \frac{(1 - \alpha)\, 2^{\alpha - 1}}{l_{B_1}} 
        \int_{-1}^{+1} \frac{[B_{w_0,xx}(\zeta)]}{(1 + \zeta)^{\alpha}} \, \mathrm{d}\zeta,
        \label{eq:B2_tilde}
    \end{equation}
        \begin{equation}
        [\tilde{B}_{w_0,yx}(\mathbf{x}, l_{A_1}, l_{B_1}, \alpha)] = 
        \frac{(1 - \alpha)\, 2^{\alpha - 1}}{l_{A_1}} 
        \int_{-1}^{+1} \frac{[B_{w_0,yx}(\zeta)]}{(1 - \zeta)^{\alpha}} \, \mathrm{d}\zeta +
        \frac{(1 - \alpha)\, 2^{\alpha - 1}}{l_{B_1}} 
        \int_{-1}^{+1} \frac{[B_{w_0,yx}(\zeta)]}{(1 + \zeta)^{\alpha}} \, \mathrm{d}\zeta,
        \label{eq:B4_tilde}
    \end{equation}
        \begin{equation}
        [\tilde{B}_{w_0,xy}(\mathbf{x}, l_{A_2}, l_{B_2}, \alpha)] = 
        \frac{(1 - \alpha)\, 2^{\alpha - 1}}{l_{A_2}} 
        \int_{-1}^{+1} \frac{[B_{w_0,xy}(\zeta)]}{(1 - \zeta)^{\alpha}} \, \mathrm{d}\zeta +
        \frac{(1 - \alpha)\, 2^{\alpha - 1}}{l_{B_2}} 
        \int_{-1}^{+1} \frac{[B_{w_0,xy}(\zeta)]}{(1 + \zeta)^{\alpha}} \, \mathrm{d}\zeta,
        \label{eq:B3_tilde}
    \end{equation}
    \begin{equation}
        [\tilde{B}_{w_0,yy}(\mathbf{x}, l_{A_2}, l_{B_2}, \alpha)] = 
        \frac{(1 - \alpha)\, 2^{\alpha - 1}}{l_{A_2}} 
        \int_{-1}^{+1} \frac{[B_{w_0,yy}(\zeta)]}{(1 - \zeta)^{\alpha}} \, \mathrm{d}\zeta +
        \frac{(1 - \alpha)\, 2^{\alpha - 1}}{l_{B_2}} 
        \int_{-1}^{+1} \frac{[B_{w_0,yy}(\zeta)]}{(1 + \zeta)^{\alpha}} \, \mathrm{d}\zeta.
        \label{eq:B5_tilde}
    \end{equation}
\end{subequations}
In this framework, $[B_{w_0,xx}(\mathbf{x})]$, $[B_{w_0,yy}(\mathbf{x})]$, and $[B_{w_0,xy}(\mathbf{x})] = [B_{w_0,yx}(\mathbf{x})]$ represent the derivative matrices associated with the local (integer-order) derivatives of the transverse displacement, as defined in Eq.~\eqref{local_derivative_approx}.  

The integrals included in these formulations are evaluated numerically using the Gauss--Jacobi quadrature technique. As a result, the fractional derivative of order $\alpha \in (0,1]$ at a spatial location $\mathbf{x}$, within the nonlocal domain $(x_j - l_{A_j},\, x_j + l_{B_j})$ along the $j$-th direction, can be expressed through the nodal values and the corresponding fractional strain--displacement matrix $[\tilde{B}]$.  

It should be noted that the governing integro-differential equations in Eqs.~\eqref{eq:inplane_gde_cpt} and \eqref{eq:transverse_GDE} incorporate both the Caputo and Riemann--Liouville definitions of fractional derivatives. The numerical formulation provided above for the Caputo derivative can be readily extended to the Riemann--Liouville derivative. However, the weak form derived from the strong form of these equations involves only Caputo-type derivatives. Therefore, the current development focuses on constructing the MLS-based approximation for 2D Caputo fractional derivatives.

\subsection{2D f-EFG model}
\label{sec:fefg_solver}

In the present section, a 2D f-EFG model is developed for nonlocal Kirchhoff plate. The formulation begins with the weak form of the governing integro-differential equations. Within this weak formulation, the fractional-order derivatives are numerically approximated in terms of nodal variables using the previously derived MLS approximants. This procedure transforms the FDEs into an equivalent system of algebraic equations for the nonlocal Kirchhoff plate.

To employ the principle of minimum potential energy for the fractional-order plate, the internal strain energy and the external work due to the applied transverse load are defined as:
\begin{subequations}
\label{eq:virtual_work}
    \begin{equation}
\delta U = 
\int_{\partial \Omega} \int_{-\frac{h}{2}}^{\frac{h}{2}} 
\left[
\tilde{\sigma}_{xx} \, \delta \tilde{\varepsilon}_{xx} +
\tilde{\sigma}_{yy} \, \delta \tilde{\varepsilon}_{yy} +
2 \tilde{\sigma}_{xy} \, \delta \tilde{\varepsilon}_{xy}
\right]
\, \mathrm{d}z\,\mathrm{d}S,
\label{eq:internal_virtual_work}
\end{equation}
\begin{equation}
\delta V =
\int_{\partial \Omega} 
\left[
f_z \, \delta w_0 
\right]
\, \mathrm{d}S.
\label{eq:external_virtual_work}
\end{equation}
\end{subequations}
The above expressions can be recast in matrix form as follows:
\begin{subequations}
    
\begin{equation}
\delta U = \int_{\partial \Omega} \int_{-\frac{h}{2}}^{\frac{h}{2}} \delta \{\tilde{\epsilon}\}^T\{\tilde{\sigma}\}~\mathrm{d}z~\mathrm{d}S,~~~\delta V=\int_{\partial \Omega} \delta \{\bm{u}_g\}^T\{f\}~\mathrm{d}S.
\end{equation}
Here, 
\begin{equation}
    \{\tilde{\epsilon}\}^T=\begin{bmatrix}
        \tilde{\epsilon}_{xx} & \tilde{\epsilon}_{yy} & \tilde{\gamma}_{xy}
    \end{bmatrix},~~~\{\tilde{\sigma}\}^T=\begin{bmatrix}
        \tilde{\sigma}_{xx} & \tilde{\sigma}_{yy} & \tilde{\sigma}_{xy}
    \end{bmatrix}
\end{equation}
are the vectors of strain and stress components at point $\mathbf{x}$, and: 
\begin{equation}
    \{\bm{u}_g\}^T=\begin{bmatrix}
        \{\bm{u}_0\}^T & \{\bm{v}_0\}^T & \{\bm{w}_0\}^T
    \end{bmatrix}
\end{equation}
represents the global vector of nodal values for axial and transverse displacement. For brevity, the explicit dependence of the field variables on the spatial coordinates within the 2D domain is omitted. The stress and strain vectors are connected through the following constitutive relation for isotropic materials:
\begin{equation}
   \label{eq: material_constt_plate}
    \{\tilde{\sigma}\}=[C]\{\tilde{\epsilon}\},~~~\text{where}~~~[C]=\begin{bmatrix}
        \frac{E}{1-\nu^2} & \frac{\nu E}{1-\nu^2} & 0\\
        \frac{\nu E}{1-\nu^2} & \frac{E}{1-\nu^2} & 0\\
        0 & 0 & \frac{E}{2(1+\nu)}.
    \end{bmatrix}
\end{equation}
\end{subequations}
The fractional-order strain components in the above expression are computed in terms of the nodal displacement values using the strain--displacement relations defined in Eq.~\eqref{eq:frac_strain_comp}. Specifically, by employing the numerical approximation of fractional derivatives based on the MLS shape functions derived earlier, the fractional-order strains can be expressed as
\begin{subequations}
\label{eq: frac_strain_mls_approx}
\begin{equation}
    \{\tilde{\epsilon}\}=[Z_2]\left[[\tilde{B}_L]+\frac{1}{2} [\tilde{B}_{N}] \right] \{\bm{u}_g\}.
\end{equation}
In the above expression, 
\begin{equation}
    [Z_2]=\begin{bmatrix}
        1 & 0 & 0 & -z & 0 & 0\\
        0 & 1 & 0 & 0 & -z & 0\\
        0 & 0 & 1 & 0 & 0 & -z
    \end{bmatrix},
\end{equation}   
and the matrices $[\tilde{B}_L]$ and $[\tilde{B}_{NL}]$ correspond to linear and nonlinear components of the fractional-order strain-displacement relations. These matrices are given by:
\begin{equation}
\label{eq_BL}
    [\tilde{B}_L]=\begin{bmatrix}
        [\tilde{B}_{u_0,x}] & [\textbf{0}] & [\textbf{0}]\\
        [\textbf{0}] & [\tilde{B}_{v_0,y}] & [\textbf{0}]\\
        [\tilde{B}_{u_0,y}] & [\tilde{B}_{v_0,x}] & [\textbf{0}]\\
        [\textbf{0}] & [\textbf{0}] & [\tilde{B}_{w_0,xx}]\\
        [\textbf{0}] & [\textbf{0}] & [\tilde{B}_{w_0,yy}]\\
        [\textbf{0}] & [\textbf{0}] & [\tilde{B}_{w_0,xy}]+[\tilde{B}_{w_0,yx}]\\
    \end{bmatrix},~~~[\tilde{B}_{N}]=\begin{bmatrix}
        [\textbf{0}] & [\textbf{0}] & D_{x}^{\alpha}w_0~[\tilde{B}_{w_0,x}]\\
        [\textbf{0}] & [\textbf{0}] & D_{y}^{\alpha}w_0~[\tilde{B}_{w_0,y}]\\
        [\textbf{0}] & [\textbf{0}] & D_{x}^{\alpha}w_0~[\tilde{B}_{w_0,y}]+D_{y}^{\alpha}w_0~[\tilde{B}_{w_0,x}]\\
        [\textbf{0}] & [\textbf{0}] & [\textbf{0}]\\
        [\textbf{0}] & [\textbf{0}] & [\textbf{0}]\\
        [\textbf{0}] & [\textbf{0}] & [\textbf{0}]\\
    \end{bmatrix}
\end{equation}
\end{subequations}
where $[\textbf{0}]$ is a $1\times n$ row matrix of zeros. The matrix $[\tilde{B}_L]$ includes terms corresponding to linear in-plane and bending strains in Eq.~\eqref{eq:frac_strain_comp}. The matrix $[\tilde{B}_{N}]$ includes nonlinear terms within the fractional-order strain-displacement relations. 

The equilibrium equations for the fractional-order Kirchhoff plate are derived from the principle of minimum potential energy \( \delta \Pi = 0 \), as previously seen in \S~\ref{subsec: gdes_frac}, for all admissible virtual displacements. We substitute the numerical approximations for fractional-order stresses and strains developed in Eqs.~\eqref{eq: material_constt_plate} and \eqref{eq: frac_strain_mls_approx} in Eq.~\eqref{eq:virtual_work} to give the following result:
\begin{equation}
    \label{eq:variation}
    \delta\Pi=\int_v\delta\{\bm{u}_g\}^T\left([\tilde{B}_L]+\frac{1}{2}[\tilde{B}_N]\right)^T[Z_2]^T[C][Z_2]\left([\tilde{B}_L]+\frac{1}{2}[\tilde{B}_N]\right)\{\bm{u}_g\}\,\mathrm{d}V -\int_{\partial \Omega}\delta\{\bm{u}_g\}^T[N_w]^T f_z \,\mathrm{d}S,
\end{equation}
where 
\begin{equation}
    [N_w]=\begin{bmatrix}
        [\textbf{0}] & [\textbf{0}] & [\phi]
    \end{bmatrix}
\end{equation}   
is the row matrix of shape functions for the transverse displacement, defined in terms of the MLS approximation functions developed in \S~\ref{subsec: mls_func}. Note that the above weak form corresponds to a numerical approximation using 2D f-EFG developed here for the geometrically nonlinear response of a fractional-order Kirchhoff plate. More clearly, the 2D MLS approximation for fractional-order derivatives reduces the integro-differential terms to algebraic expressions. 

Applying the fundamental lemma of variational calculus ($\delta \Pi=0$), we get the following nonlinear algebraic equations for the equilibrium of a fractional-order Kirchhoff plate.
\begin{equation}
    \label{eq:algebric}
    [K]\{\bm{u}_g\}=\{F\}.
\end{equation}
In the above equation, $[K]$ represents the global stiffness matrix of the plate and ${F}$ is the transverse mechanical load vector applied to the plate. The system matrices in the above equation are given as follows:
\begin{subequations}
\label{eq:matrix}
\begin{equation}
    \label{eq:stiffness}
    [K]=\int_{\partial \Omega}\left[[\tilde{B}_L]^T[D_t][\tilde{B}_L]+\frac{1}{2}[\tilde{B}_N]^T[D_t][\tilde{B}_L]+\frac{1}{2}[\tilde{B}_L]^T[D_t][\tilde{B}_N]+\frac{1}{4}[\tilde{B}_N]^T[D_t][\tilde{B}_N]\right]\,\mathrm{d}S,
\end{equation}
\begin{equation}
    \{F\}=\int_{\partial \Omega} [N]^T\{f\}~\mathrm{d}S.
\end{equation}
 Here,
     \begin{equation}
         \label{eq:dt}
         [D_t]=\int_{-h/2}^{h/2}[Z_2]^T[C][Z_2]\,\mathrm{d}z.
     \end{equation}
\end{subequations}
The above system matrices present the 2D f-EFG numerical model for nonlinear (spatial) FDEs. while, the FDEs are chosen here to model the nonlocal behavior of Kirchhoff plates. The developed numerical method is general enough to be applied to other FDEs beyond elasticity. These system matrices are computed numerically using the Gauss-Legendre quadrature, whereas the Gauss-Jacobi quadrature was used earlier for the fractional-order derivatives. A background mesh is constructed for this integration, typically aligned with the nodal distribution used in the MLS approximation. 
Note that the expression for $[B_N]$ depends on the deformation state of the plate, rendering the system matrices nonlinear. The nonlinear algebraic system is solved using an iterative Newton-Raphson method, with incremental loading applied to improve convergence. For a detailed description of the incremental iterative Newton-Raphson procedure, the reader is referred to~\cite{sidhardh2020geometrically}.

Furthermore, it should be emphasized that the approximation functions constructed by MLS interpolation do not satisfy the Kronecker delta property \cite{liu2005introduction,belytschko1994element}. As a result, the exact imposition of essential boundary conditions is not directly possible within the f-EFG framework. To address this, a point collocation technique is adopted here to impose the essential boundary conditions.\\

\noindent
\textit{\underline{\textbf{Linear 2D f-EFG}}}: The nonlinear algebraic system presented earlier characterizes the finite deformation response of the nonlocal plate. The strain-displacement relations in Eq.~\eqref{eq: frac_strain_mls_approx} are modified to retain only the linear terms, that is, by omitting $[B_N]$, and applying the corresponding changes to the weak form in Eq.~\eqref{eq:virtual_work}, the 2D f-EFG formulation in Eq.~\eqref{eq:algebric} reduces to a linear algebraic system representing the plate under infinitesimal deformations. A detailed derivation of the f-EFG system matrices for this linear elastic response of a nonlocal plate is not repeated here for the sake of brevity.

\section{Results and Discussion}
\label{sec: results}
In the previous section, a 2D f-EFG solver was developed for linear and nonlinear FDEs. In this section, the solver is employed to conduct a numerical study of the linear and geometrically nonlinear response of a nonlocal Kirchhoff plate. The investigation begins with an assessment of the convergence of the 2D f-EFG solver, followed by a validation with results reported in the literature. Additional numerical examples are then presented to demonstrate the accuracy and performance of the proposed f-EFG formulation.
  
In the present analysis, the nonlocal plate is subjected to a UDTL, with its magnitude $f_z(x,y)=q_0$;~~$q_0$ (in Pa) is specified whenever necessary. The formulation developed in the previous section is based on the classical plate theory, also known as the Kirchhoff plate displacement theory, which assumes a sufficiently large aspect ratio ($a/h$ or $b/h$ $\geq 50$). Therefore, the geometry of the rectangular plate is chosen to be: length $a = 1~\mathrm{m}$ and width $b = 1~\mathrm{m}$, while the thickness is  $h = 0.02~\mathrm{m}$ \cite{patnaik2021fractional}. The plate is assumed to be composed of an isotropic material characterized by a Young’s modulus of $E = 1.09~\mathrm{MPa}$ and a Poisson's ratio of $\nu=0.3$\cite{reddy2015introduction}.The fractional-order constitutive parameters $\alpha$, $l_A$, and $l_B$, which capture long-range interactions within the nonlocal plate, are varied and specified as needed.

Recall that the $l_A$ and $l_B$ corresponding to the nonlocal horizon of influence depend on the position of the point in the 2D domain. However, for the isotropic material considered in this study, these length scales are assumed to be equal and constant at points sufficiently far from the boundaries, denoted $l_{{A}{\Box}} = l_{{B}{\Box}} = h_{l}, \quad \Box \in {x,y}$.
This constant $h_l$ represents the intrinsic/characteristic material length scale; assumed to be identical along the $x-$ and $y-$ directions for the isotropic material. However, for points close to the boundaries of the 2D domain, the nonlocal lengths along the $x-$ and $y-$ directions are appropriately truncated \cite{patnaik2021fractional}. It is important to emphasize that the use of position-dependent length scale parameters is motivated by the need to maintain frame invariance and ensure kernel completeness in the presence of asymmetric horizons, material boundaries, and interfaces \cite{sumelka2014thermoelasticity,patnaik2020ritz}. This consideration is especially critical in the study of fracture and damage in solids, where internal material boundaries can disrupt long-range interactions across the domain. For complete details on the position-dependent length scales, please refer to \cite{patnaik2021fractional}. However, the numerical framework developed here is flexible and can be extended to accommodate alternative specifications for the length scales $l_A$ and $l_B$ as discussed in \cite{sumelka2017fractional, sumelka2014fractional}. 

The transverse displacements of the plate, evaluated under the assumption of infinitesimal deformations using the linear 2D f-EFG model, are non-dimensionalized as follows to facilitate comparison with the existing literature\cite{reddy2015introduction}:
\begin{equation}
    \label{eq:normalization}
    \overline{w} (\mathbf{x})=w_0\left(\mathbf{x}\right)\left[\frac{100Eh^3}{q_0a^4}\right]
\end{equation}
Furthermore, the transverse displacement and the applied load for the nonlinear response are reported as \cite{reddy2015introduction}:
\begin{equation}
    \label{eq:normalization_1}
    \bar{w}(\mathbf{x})=\frac{w_0\left(\mathbf{x}\right)}{h},~~~~~~~~~~~~\bar{P}=\frac{q_0a^4}{Eh^3},
\end{equation}
 in the convergence study. All symbols in the above expressions are defined previously, and the transverse displacement is evaluated at a point $\mathbf{x}$ of the 2D domain (mid-plane of the plate).

\subsection{Convergence}
\label{subsec:Convergence}
The convergence of the 2D f-EFG solver developed here is examined by comparing the transverse displacement of the rectangular plate for different numbers of nodes. More specifically, the nonlocal rectangular plate is discretized by an increasing number of nodes, while ensuring sufficient number within the support domain of the MLS approximation. The latter condition is essential to ensure that $[A(\mathbf{x})]$ in Eq.~\eqref{eq:a_solution} is nonsingular. The nonlocal plate is characterized by fractional order $\alpha = 0.8$ and length scale $h_l = 0.5a$. 

We note that unless otherwise specified, the nodes for the 2D f-EFG model are uniformly distributed across the 2D domain. More clearly, the nodes are distributed across the 2D rectangular domain so that the distance between any two neighboring nodes along the $x-$ and $y-$ directions is given by $l_{ex}$ and $l_{ey}$. We repeat the numerical simulations for different choices of these lengths to check for convergence. 

\begin{figure}[H]
    \centering
    \includegraphics[width=0.5\linewidth]{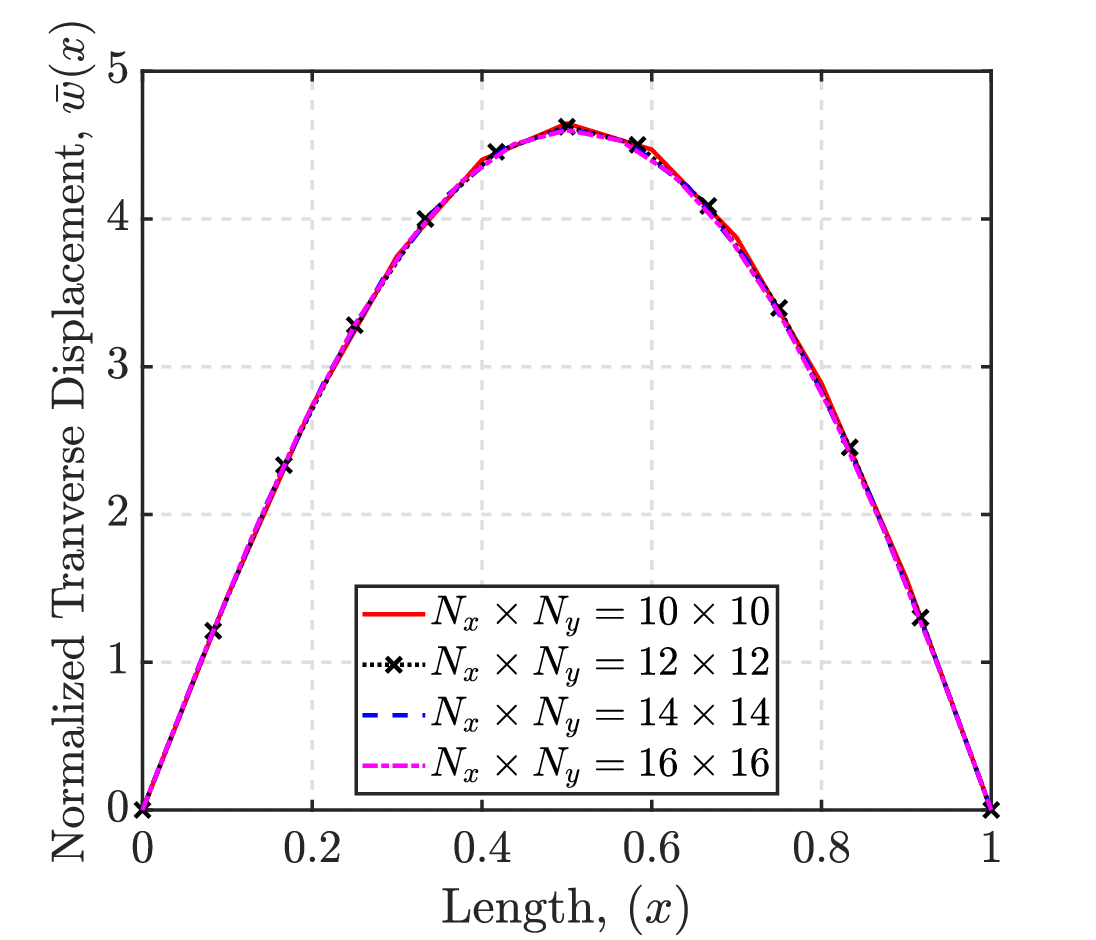}
    \caption{Transverse displacement for the linear (infinitesimal) response of nonlocal rectangular plate compared for different node distributions.}
    \label{fig:Convergence}
\end{figure}

Recall that the f-EFG solver includes an additional integration step to account for the differ-integral definition of fractional-order derivatives. This step is due to a nonlocal definition of the matrix $[\tilde{B}]$ (see Eq.~\eqref{b2 b3 tildes}), unlike standard EFG solvers for integer-order differential equations. In the 2D f-EFG formulation presented here, the integration is carried out independently along the $x$- and $y$-directions. A sufficient density of nodes is required within the nonlocal horizon along both directions to achieve an accurate numerical evaluation of the fractional-order derivatives. To quantify this, the \textit{dynamic rate of convergence} is defined as \cite{patnaik2021fractional,rajan2024element}.
\begin{equation}
\label{eq: dynamic_rate}
\mathcal{N}_x \times \mathcal{N}_y \equiv 
\left(\frac{h_l}{l_{ex}} \times \frac{h_l}{l_{ey}}\right).
\end{equation}
This choice ensures a sufficient number of nodes within the nonlocal influence horizon. Figure~\ref{fig:Convergence} illustrates a comparison of the linear response of the fractional-order Kirchhoff plate. More clearly, the normalized transverse displacement $\overline{w}(x,b/2)$ predicted by the 2D f-EFG model is compared for different node distributions to assess the rate of convergence. The comparison indicates that the 2D f-EFG model converges, with successive refinements in the nodal discretization, leading to variations less than $2\%$. Moreover, it is found that convergence is achieved for a nodal configuration of $\mathcal{N}_x \times \mathcal{N}_y = 12 \times 12$ for the fractional-order Kirchhoff plates.

To examine the convergence of the nonlinear f-EFG model, the rectangular plate is analyzed for a UDTL of $\bar{P} = 250$ (see Eq.~\eqref{eq:normalization_1} for normalization definition). The transverse displacement at mid-point is calculated for different fractional-order constitutive parameters ($\alpha$ and $h_l$) and for various $\mathcal{N}_x \times \mathcal{N}_y$. It is noted that, for the nonlinear root-finding iterations, the tolerance in the $\mathbb{L}^2$ norm of incremental displacement is prescribed as $10^{-3} \times h$. The corresponding results are compared in Table~\ref{tab:convergence}.

\begin{table}[H]
    \centering
    \begin{tabular}{c|c|c c c c}
    \hline\hline
    & \multirow{2}*{$\mathcal{N}_x \times \mathcal{N}_y$} & \multicolumn{4}{c}{Normalized transverse displacement, $\bar{w} (a/2,b/2)$}\\
        &  & $\alpha=1.0$ &  $\alpha=0.9$ &  $\alpha=0.8$ & $\alpha=0.7$ \\
        \hline\hline
        \multirow{4}*{$h_l=0.5a$} & $12\times12$ & 1.834 & 1.950  & 2.068 & 2.199 \\
        & $14\times14$ & 1.835 & 1.946 & 2.042 & 2.267 \\ 
        & $16\times16$ & 1.832 & 1.945 & 2.046 & 2.268 \\
        & $18\times18$ & 1.826 & 1.949 & 2.062 & 2.299\\ 
         \hline
                \multirow{4}*{$h_l=0.6a$} & $12\times12$ & 1.834 & 2.004 & 2.233 &  2.537 \\
        & $14\times14$ & 1.835 & 2.004 & 2.218 & 2.520\\ 
        & $16\times16$ & 1.832 & 2.005 & 2.232 & 2.547\\
        & $18\times18$ & 1.826 & 1.997 & 2.212 &  2.530\\ 
         \hline\hline
    \end{tabular}
    \caption{Normalized transverse displacement for the nonlinear response of nonlocal rectangular plate compared for different node distributions.}
    \label{tab:convergence}
\end{table}
The nonlinear f-EFG models demonstrate convergence beyond $\mathcal{N}_x \times \mathcal{N}_y= 12\times12$, as the results change by less than $2\%$. This observation holds for all choices of fractional-order constitutive parameters considered in the study. Consequently, a uniform distribution of $\mathcal{N}_x \times \mathcal{N}_y= 12\times12$ nodes is adopted for the validation studies (linear and nonlinear) on 2D rectangular plate domain, except in Section~\ref{subsec: nonuniform}, where a non-uniform distribution is explored.

\subsection{Validation}
In this section, the linear and geometrically nonlinear responses of nonlocal plates obtained using the 2D f-EFG solver are compared with the finite element results reported in \cite{patnaik2021fractional,patnaik2020geometrically}. This comparison highlights the influence of the fractional-order constitutive parameters, $\alpha$ and $h_l$.\\

\noindent
\textbf{Validation \#1}: The linear response of the rectangular plate subjected to an UDTL is solved using the 2D f-EFG numerical solver. The transverse displacements at the midpoint $\mathbf{x}(a/2,b/2)$ in the mid-plane are non-dimensionalized following Eq.~\eqref{eq:normalization} and compared with \cite{patnaik2021fractional}. This comparison, carried out for various combinations of the fractional-order parameters $\alpha$ and $h_l$, is shown in Table~\ref{tab:linear_validation}. The excellent agreement between the predictions of the 2D f-EFG solver and the reference data, with an error $<1\%$, confirms the validation of the 2D f-EFG solver for the linear elastic response of the nonlocal plate. Furthermore, this excellent agreement points to the validation of the 2D f-EFG solver for a numerical solution of 2D linear FDEs.
\begin{table}[H]
    \centering
    \begin{tabular}{c|c|c c c c}
    \hline\hline
     & &\multicolumn{4}{c}{Normalized transverse displacement, $\overline{w}(a/2,b/2)$}\\
                    &   & $\alpha=1.0$ &   $\alpha=0.9$ &   $\alpha=0.8$ &$\alpha=0.7$ \\
        \hline\hline
        \multirow{2}*{$h_l=0.5a$} & \cite{patnaik2021fractional}& 4.5701 & 4.9480 & 5.4094 & 6.0180 \\
        & f-EFG & 4.5844 & 4.9462 & 5.4110 & 6.0426\\
         \hline
          \multirow{2}*{$h_l=0.3a$} & \cite{patnaik2021fractional}& 4.5701 & 4.7068 & 4.8249 & 4.8768 \\
        & f-EFG & 4.5844 & 4.7042 & 4.8249 & 4.8768 \\
        \hline 
        \hline
    \end{tabular}
    \caption{Normalized transverse displacement of rectangular plate predicted here is compared with f-FEM~\cite{patnaik2021fractional} for different nonlocal constitutive parameters. Infinitesimal deformations are considered here.}
    \label{tab:linear_validation}
\end{table}

\noindent
\textbf{Validation \#2}: The geometrically nonlinear response for a rectangular plate under an applied UDTL of $q_0=10$~Pa is obtained using the 2D f-EFG solver developed in \S~\ref{sec:fefg_solver}. The nonlocal response is computed for different values of the fractional-order constitutive parameters $\alpha$ and $h_l$, and the results are summarized in Table~\ref{tab:non_lin_validation}. In addition, these results are compared with those reported in the literature using a nonlinear FE solver for FDEs \cite{patnaik2020geometrically}. The strong agreement between the two result sets confirms the accuracy and efficacy of the 2D f-EFG numerical method developed here to solve nonlinear fractional-order boundary value problems.
\begin{table}[H]
    \centering
    \begin{tabular}{c|c|c c c c}
    \hline\hline
    & &\multicolumn{4}{c}{Transverse displacement, $w(a/2,b/2)$ (in m)}\\
                    &   & $\alpha=1.0$ &   $\alpha=0.9$ &   $\alpha=0.8$ &   $\alpha=0.7$ \\
        \hline\hline
        \multirow{2}*{$h_l=0.6a$} & \cite{patnaik2020geometrically}& 0.1029 & 0.1135 & 0.1262 & 0.1423 \\
        & f-EFG &  0.1052 & 0.1145 & 0.1260 & 0.1413\\
         \hline
          \multirow{2}*{$h_l=0.5a$} & \cite{patnaik2020geometrically}&  0.1029 & 0.1096 & 0.1170 & 0.1269  \\
        & f-EFG &  0.1052 & 0.1103 & 0.1154 & 0.1221 \\
        \hline\hline
    \end{tabular}
    \caption{Normalized transverse displacement of rectangular plate predicted here is compared with f-FEM~\cite{patnaik2020geometrically} for different nonlocal constitutive parameters. Finite deformations are considered here.}
    \label{tab:non_lin_validation}
\end{table}

\subsection{Non-uniform grid}
\label{subsec: nonuniform}
We have already established the efficacy and accuracy of the 2D f-EFG method using the MLS approximation for nonlocal plates with uniformly distributed nodes. Unlike the finite element method, which depends on mesh connectivity, the MLS approximation relies solely on the distribution of nodes within the domain. This is particularly relevant for problems such as adaptive domain geometry, crack propagation \cite{belytschko1994fracture}, and damage \cite{askes2000dispersion} that require extensive remeshing following FEM. In contrast, meshfree methods are particularly well suited for such applications due to the absence of predefined element connectivity \cite{belytschko1999smoothing}.  

\begin{figure}[H]
    \centering
    \includegraphics[width=1\linewidth]{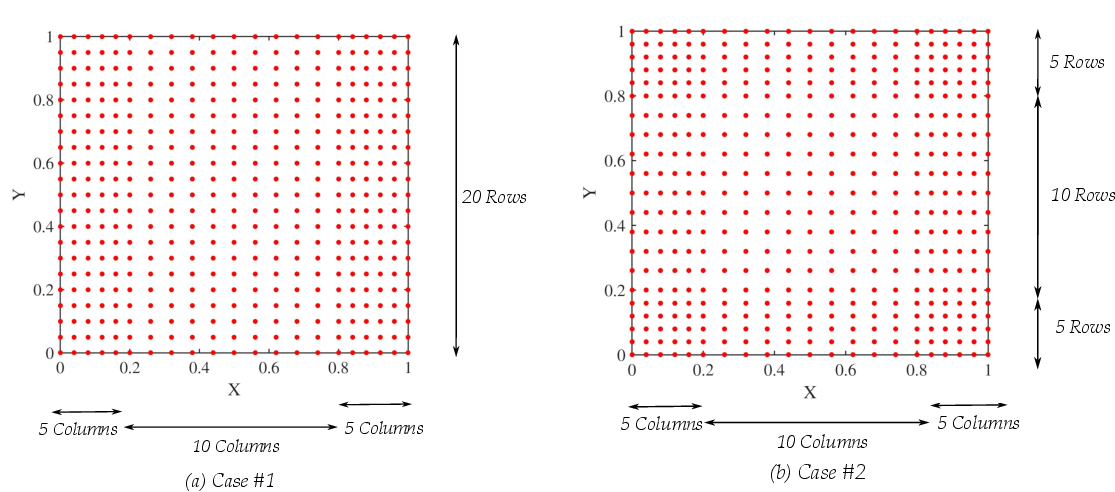}
    \caption{A schematic representation of the non-uniform grid in a rectangular plate.}
    \label{fig:Nonuniform_mesh}
\end{figure}
\begin{figure} [H]
    \centering
    \begin{subfigure}{.5\textwidth}
        \centering
        \includegraphics[width=1\linewidth]{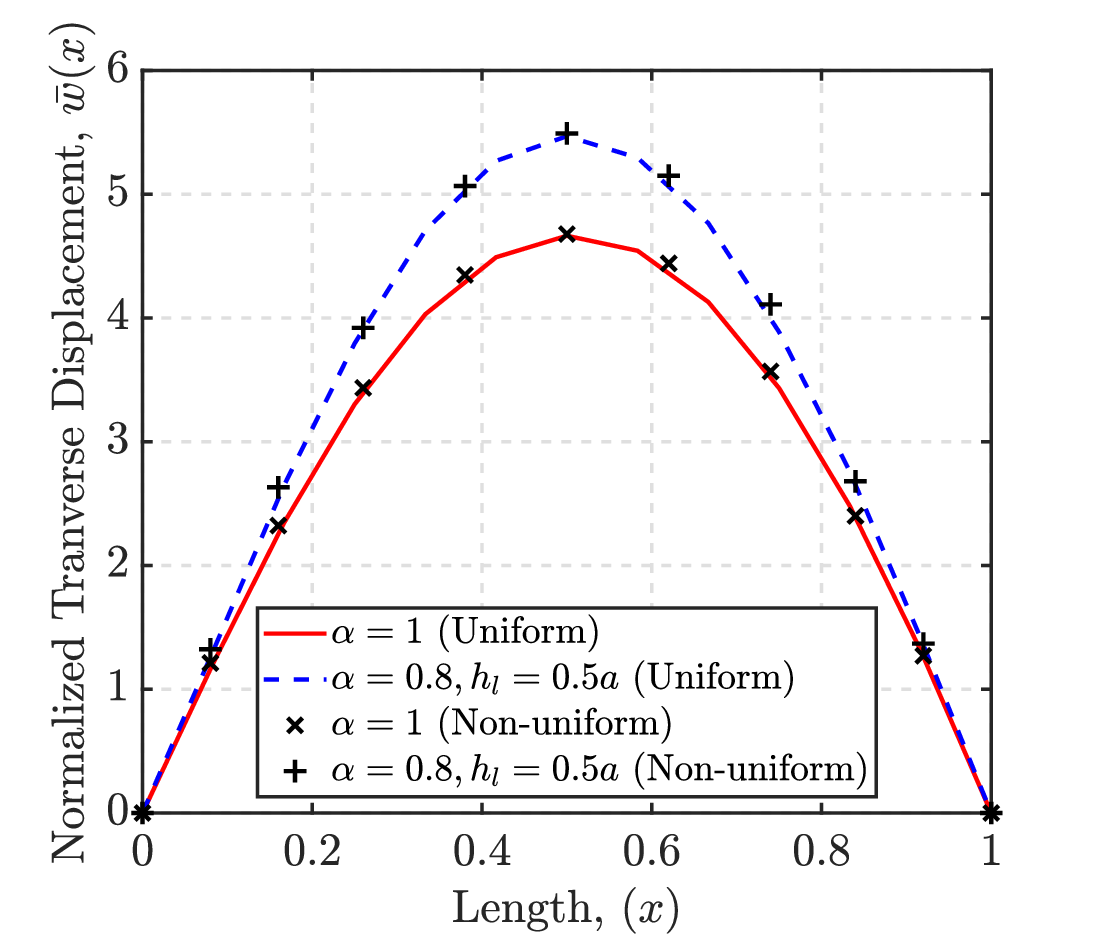}
        \label{fig:simply_validation}
    \end{subfigure}%
    \begin{subfigure}{.5\textwidth}
        \centering
        \includegraphics[width=1\linewidth]{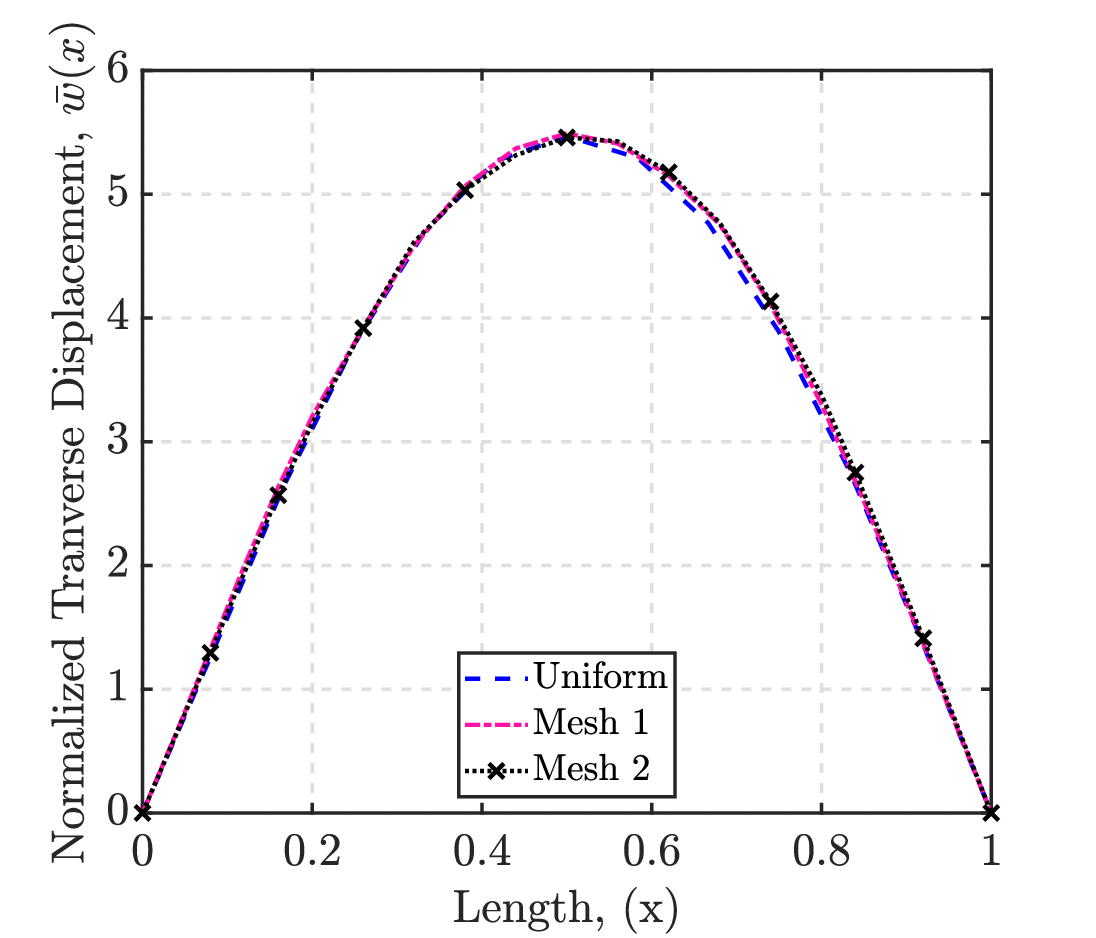}
        \label{fig:simply_variation_non_uniform_mesh}
    \end{subfigure}
    \caption{Normalized transverse displacement along the length of rectangular plate at $y=b/2$ under a UDTL. The linear elastic response is compared for both uniform and non-uniform node distributions.}
    \label{fig:nonuniform}
\end{figure}

In this study, we investigate the efficacy of the 2D f-EFG numerical solver for non-uniform node distributions, such as those illustrated in Fig.~\ref{fig:Nonuniform_mesh}. Specifically, the nonlocal Kirchhoff rectangular plate is solved using the 2D f-EFG model developed above, with the two different choices of non-uniform node distributions, and the results are compared against those obtained for a uniform node distribution. In Case \#1, the nodes are distributed uniformly along the $y$-direction but not uniformly along the $x$-direction. This configuration is particularly relevant for problems involving material heterogeneity (functionally graded) or notches, where localized refinement is needed only in one direction \cite{davoudi2018comparison}. In Case \#2, the nodes are distributed non-uniformly along both directions. The corresponding results for non-uniform and uniform node distributions are compared in Fig.~\ref{fig:nonuniform}.

The MLS approximation, together with the f-EFG formulation, is constructed for two non-uniform mesh distributions over the 2D plate domain. In each case, a non-uniform support domain ($d_{m}$) is adopted in the MLS approximation to account for the varying node density within the domain. The remaining steps in developing the numerical solver for non-uniform node distribution remain the same. The f-EFG model is limited to the linear elastic response of a rectangular plate subjected to a UDTL. In Fig.~\ref{fig:nonuniform}a, the transverse displacement field is computed for two cases: (i) the local elastic model with $\alpha = 1$ and (ii) a fractional-order elastic model with $\alpha = 0.8$ and length scale $h_l = 0.5a$. The displacement along the length predicted by the non-uniform mesh depicted as Case~\#1 in Fig.~\ref{fig:Nonuniform_mesh}a is compared with that of a uniform $12 \times 12$ distribution of nodes; ($a/l_{ex}=b/{l_{ey}}=12$). The transverse displacements are non-dimensionalized according to Eq.~\eqref{eq:normalization} and compared in Fig.~\ref{fig:nonuniform}. The comparison shows excellent agreement, with discrepancies below 1$\%$, between results for non-uniform and uniform distribution of nodes for integer- and fractional-order constitutive parameters. This observation highlights the robustness and accuracy of the f-EFG solver even for non-uniform node distributions in a 2D plate.
To further evaluate the impact of non-uniform mesh distributions, the linear elastic response of the rectangular plate is analyzed using the two mesh distributions shown in Fig.~\ref{fig:Nonuniform_mesh}. The f-EFG model is solved for uniform and non-uniform node distributions considering the fractional-order constitutive parameters as: $\alpha = 0.8$ and $h_l = 0.5a$. The resulting transverse displacement fields are compared in Fig.~\ref{fig:nonuniform}, and demonstrate excellent consistency across all cases. However, it is important to note that a higher total number of nodes may be required to achieve convergence in the EFG framework following non-uniform distribution of nodes \cite{krysl1995analysis,rajan2024element}.

The above observation of consistent results for the 2D f-EFG solver using non-uniform distribution of nodes is significant. 
This is due to the challenges of using a non-uniform distribution of nodes for the MLS approximation \cite{liu2005introduction,li2016stability}. Some of the challenges of using non-uniform distributions are the near-singular and ill-conditioned $[A]$ matrix (see Eq.~\eqref{eq:A_H_definitions}), a simultaneous increase in computational cost, and instability of the approximations~\cite{li2016stability}. These issues reported in the literature for integer-order derivatives are only exaggerated for fractional-order derivatives. However, as discussed previously, such a non-uniform distribution of nodes is relevant for practical problems. Therefore, the consistent performance of the 2D f-EFG solver with non-uniform distribution of nodes demonstrates the efficacy of the 2D f-EFG solver for fractional-order boundary value problems. 

\subsection{Circular plate}
Firstly, we note that the bending response of circular plates following fractional-order constitutive relations for nonlocal elasticity has not been explored in the literature. This is owing to the challenges in the evaluation of the fractional-order derivatives in an irregular domain following a mesh-based numerical solver. In this study, we address these challenges by applying the meshfree f-EFG solver to numerically solve fractional-order boundary value problems in complex, irregular geometries. Moreover, as mentioned previously, this advantage of the meshfree solver for irregular geometries is demonstrated for 2D domains, and was not apparent using the 1D f-EFG solver developed previously\cite{rajan2024element}.

In this study, a simply supported circular plate subjected to a UDTL with a radius of $a=1~\mathrm{m}$ and thickness $h=0.02~\mathrm{m}$ is considered. A schematic illustration of the circular plate is provided in Fig.~\ref{fig:SSSS_CIRC_plate}. The plate is assumed to be composed of isotropic material with constants as considered in above studies. The steps involved in the derivation of the strong and weak form of the governing equations for the fractional-order circular Kirchhoff plates are an extension of those provided in \S~\ref{subsec: gdes_frac} for the rectangular plate; therefore, not repeated here in the interest of brevity. The f-EFG numerical model for the circular plate is also given by Eq.~\eqref{eq:algebric}, where the system matrices are evaluated over the 2D circular domain. In the derivation of the MLS approximation functions, a non-uniform distribution of nodes is considered along the Cartesian $x-$ and $y-$directions. The choice of nodes along the $x-$ and $y-$ directions is considered to ensure compatibility with the definition of fractional-order derivatives along the $x-$ and $y-$ directions. The uniform distribution of nodes within a circular domain either leads to a poor approximation of the circular boundary geometry or requires a very large number of nodes to achieve an accurate representation. Therefore, a gradually increasing discretization is considered while approaching the circular boundaries to obtain an improved geometric approximation of the circular boundary. The efficacy of the f-EFG solver for non-uniform node distribution is already established above.

We begin with the linear response of the circular plate solved using the 2D f-EFG solver. A UDTL of $q_0=1$~Pa is applied throughout the surface of the plate and the corresponding transverse displacement at the mid-plane is evaluated in the center of the plate. As validation of the 2D f-EFG model developed here for cicular geometry, the maximum transverse displacement of the circular plate for $\alpha=1.0$ is compared with the analytical solution for the local elastic plate given in \cite{reddy2006theory}. The analytical expression for the maximum transverse displacement at the mid-plane in the center of the plate ($\mathbf{x}(0,0)$) is given by \cite{reddy2006theory}:
\begin{equation}
\label{eq:ana_circular}
w_{\text{max}} = \frac{3 q_0 a^4 (1-\nu^2)}{16 E h^3} \left( \frac{5 + \nu}{1 + \nu} \right).
\end{equation}

\begin{figure}[H]
    \centering
    \includegraphics[width=0.6\textwidth]{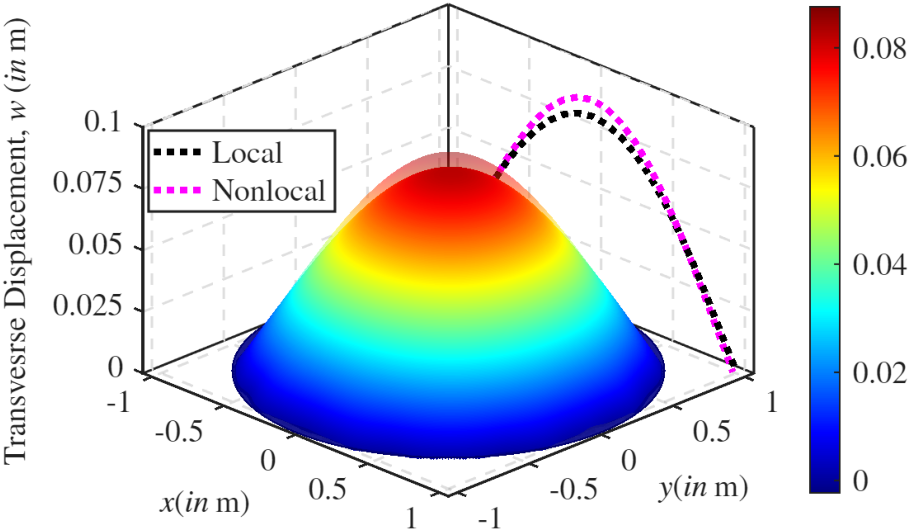}
    \caption{
    Comparison of the transverse displacement profiles for circular plate subject to UDTL considering local elastic (filled color) and nonlocal elastic response (partial opacity).
    }
    \label{fig:circ_plate_contour22}
\end{figure}

\begin{table}[H]
    \centering
    \begin{tabular}{c|c c c c c}
        \hline
        & \multicolumn{5}{c}{Maximum transverse displacement, $w_0(0,0)$} (in $\mathrm{m})$ \\
        $h_l$ & Analytical & $\alpha=1.0$ & $\alpha=0.9$ & $\alpha=0.8$ & $\alpha=0.7$ \\
        \hline\hline
        $0.8a$      & 0.0797 & 0.0809 & 0.0819 & 0.0850 & 0.0942 \\
        $a$      & 0.0797 & 0.0809 & 0.0830 & 0.0876 & 0.1011 \\
        $1.2a$   & 0.0797 & 0.0809 & 0.0841 & 0.0916 & 0.1074 \\
        \hline
    \end{tabular}
    \caption{Transverse displacement at the center of the circular plate corresponding to linear elastic response for different nonlocal constitutive parameters. Results from the 2D f-EFG model and the analytical solution\cite{reddy2006theory} are reported.}
    \label{tab:linear_validation_circle}
\end{table}
The excellent agreement between analytical and numerical results for the local elastic circular plate with error $<2\%$ demonstrates the efficacy of the 2D f-EFG solver for complex and irregular geometries. For a clearer illustration of the nonlocal effect on linear elastic response, the transverse displacement contours and mid-plane displacement graphs of the circular plate are presented in Fig.~\ref{fig:circ_plate_contour22}. The lighter opacity contour represents the transverse displacement of the nonlocal plate with $\alpha=0.8$ and $h_l=a$, while the other contour corresponds to the local response. Based on the results presented in Table~\ref{tab:linear_validation_circle} and Fig.~\ref{fig:circ_plate_contour22}, the transverse displacement increases with an increase in the degree of nonlocal interactions: a reduction in the fractional order $\alpha$ or an increase in the horizon of nonlocal influence $h_l$. The consistent softening behavior of the circular plate, observed with a decrease in $\alpha$ or an increase in $h_l$, aligns well with the findings reported for rectangular plates in the literature~\cite{patnaik2021fractional}. This observation attests to the consistency of the fractional-order framework for nonlocal elasticity, hitherto unexplored for irregular geometries.
Furthermore, the 2D f-EFG solver is extended for the geometrically nonlinear response of the circular plate. The distribution of nodes considered in the linear study is retained here. As in the case of a rectangular plate, the system matrices in this case are nonlinear, and an iterative Newton-Raphson scheme with incremental loading is employed.  A UDTL of $q_0=1$~Pa is applied throughout the surface of the plate and the corresponding transverse displacement at the mid-plane is evaluated. The transverse displacement at mid-point is reported in Table~\ref{tab:nonlin_circular} for different choices of the fractional-order constitutive parameters. 
\begin{table}[H]
    \centering
    \begin{tabular}{c| c c c c}
        \hline
        & \multicolumn{4}{c}{Maximum transverse displacement, ${w}(0,0)$ (in $\mathrm{m}$)} \\
        $h_l$  & $\alpha=1.0$ & $\alpha=0.9$ & $\alpha=0.8$ & $\alpha=0.7$ \\
        \hline\hline
        $a$       & 0.02464 & 0.02692 & 0.02939 & 0.03246  \\
        $1.2a$    & 0.02464 & 0.02760 & 0.03121 & 0.03566 \\ 
        \hline
    \end{tabular}
    \caption{Transverse displacement at the center of the circular plate corresponding to geometrically nonlinear response varying nonlocal constitutive parameters.}
    \label{tab:nonlin_circular}
\end{table}

Similarly to the linear elastic response, we note a consistent softening in the geometrically nonlinear and elastic response of the fractional-order circular plate with reducing $\alpha$ and increasing $h_l$.

\subsection{Computational cost}
This section analyzes the computational cost associated with the nonlocal 2D f-EFG method. Compared to local formulations, nonlocal simulations require significantly higher computational effort, primarily due to the increased number of operations required for stiffness matrix assembly. 
To assess this computational overhead, the number of Floating Point Operations (FLOPs) required for the evaluation of the stiffness matrix is estimated.
To provide a clear assessment, the discussion is restricted to the linear case, thereby avoiding complexities associated with nonlinear analyses, such as iterative convergence procedures or incremental loading schemes. The analysis focuses specifically on estimating the FLOPs required to construct the linear consolidated stiffness matrix, as defined in Eq.~\eqref{eq:stiff}. Only the computational cost of stiffness matrix evaluation is considered, since the procedures for computing the force vector in nonlocal formulation are analogous to those in classical systems and require a comparable number of FLOPs. In this calculation, we consider $n$ uniformly distributed nodes in the 2D domain, and the domain is discretized into $N_e$ elements using a background mesh, used only for numerical integration and not for interpolation. Within each element, Gauss-Legendre quadrature is employed using $N_{GP}$ integration points. This procedure facilitates the construction of the global stiffness matrix for the nonlocal plate formulation. 

The numerical evaluation of the stiffness matrix starts with the evaluation of the strain-displacement matrices with fractional-order derivatives in Eq.~\eqref{eq_BL}, and followed by their assembly to derive the stiffness matrix in Eq.~\eqref{eq:stiffness}.
We start by estimating the FLOPs involved in evaluation of the nonlocal strain-displacement matrices. This matrix is recast in an integral-form defined at a Gauss-Legendre quadrature point $x^j$ following Eq.~\eqref{eq:frac_derivatives} as:
\begin{equation}
    [\tilde{B}_{L}(x^j,l_A,l_B,\alpha)]=\frac{(1-\alpha)2^{\alpha-1}}{l_A}\int_{-1}^{1} \frac{[B_L(\zeta)]}{(1-\zeta)^\alpha}\mathrm{d}\zeta+\frac{(1-\alpha)2^{\alpha-1}}{l_B}\int_{-1}^{1} \frac{[B_L(\zeta)]}{(1+\zeta)^\alpha}\mathrm{d}\zeta.
\end{equation}
The above computation is repeated for all $N_{GP}$ quadrature points in the 2D domain. 
The weakly-singular integrals (at terminals) in the above equation are evaluated using Gauss-Jacobi quadrature as follows:
\begin{equation}
\label{eq:comp_cost1}
    [\tilde{B}_{L}(x^j,l_A,l_B,\alpha)]=\frac{(1-\alpha)2^{\alpha-1}}{l_A}\sum\limits_{i=1}^{N_{GJP}} \underbrace{\widehat{W}_i {[B_L(\zeta_i)]}}_{6\times G_{DOF}}+\frac{(1-\alpha)2^{\alpha-1}}{l_B}\sum\limits_{i=1}^{N_{GJP}}\underbrace{\widehat{W}_i {[B_L(\zeta_i)]}}_{6\times G_{DOF}}.
\end{equation}
The cost of matrix operations is indicated in the under-braces. In the above result, the coordinates and weight of the Gauss-Jacobi quadrature are $\zeta_i$ and $\widehat{W}_i$, for $i = 1, 2, \dots, N_{GJP}$. The total number of Gauss-Jacobi quadrature points is denoted by $N_{GJP}$. The number of FLOPs required to evaluate the integrand at a single Gauss-Jacobi quadrature point is simply $2 \times 6 \times G_{DOF}$;  derived from two steps for left- and right-fractional derivative, each requiring approximately $6 \times G_{DOF}$ FLOPs.
 Recall that $G_{DOF}=3\times n$ is the number of global degrees of freedom of the plate. Hence, performing this operation across all $N_{GJP}$ Gauss-Jacobi points results in a cost of $N_{GJP} \times 2\times 6 \times G_{DOF}$ FLOPs. This step must be repeated at $N_{\text{GP}}$ Gauss-Legendre quadrature points for each of the $N_e$ elements in the domain, As a result, the total computational cost for the evaluation of $[\tilde{B}_L]$ is $N_e\times N_{GP}\times N_{GJP}\times (2 \times 6 \times G_{DOF})$. It is important to note that this expression only accounts for the numerical evaluation of the strain-displacement derivative matrices Eq.~\eqref{eq_BL} and excludes the additional cost associated with computing the nonlinear matrices. Assuming that the mesh is sufficiently refined, such that the number of elements satisfies $\mathcal{O}(G_{DOF})\propto \mathcal{O}(N_e)$, the total computational cost in terms of FLOPs scales linearly with the global degrees of freedom. Hence, the overall computational cost of this step can be expressed as $\mathcal{O}(G_{DOF}^2)$.

The nonlocal (linear) stiffness matrix following the steps given in Eq.~\eqref{eq:stiffness} is evaluated as:
\begin{equation}
    \label{eq:stiff}[K]=\sum\limits_{i=1}^{N_{e}}~J_{i}~\sum\limits_{j=1}^{N_{GP}} \underbrace{{w}_{j}   \underbrace{[\tilde{B}_{L}(\mathbf{x}^{j})]^T \underbrace{[D_t][\tilde{B}_{L}(\mathbf{x}^{j})]}_{6\times G_{DOF}\times(2\times 6-1)}}_{(G_{DOF})^2\times(2\times 6 -1)}}_{(G_{DOF})^2},  
\end{equation}
$w_j$, $j = 1, 2, \dots, N_{GP}$ in the above equation is the quadrature weight of the $j$-th Gauss-Legendre quadrature point, and Jacobian $J_i$ accounts for the transformation from the natural coordinate system to the Cartesian coordinate system for the $i$-th 2D background element. 

The global stiffness matrix is evaluated by adding the integrand contributions at every quadrature point for all elements. The contribution of each Gauss-Legendre quadrature point ($\mathbf{x}^j$ in the $i$-th element) to nonlocal stiffness is given by $(G_{DOF})^2+(G_{DOF})^2\times(2\times 6 -1)+6\times G_{DOF}\times(2\times 6-1)$ FLOPs. Thus, the total cost for the estimation of the nonlocal stiffness of the 2D domain is $N_e\times N_{GP}\times (G_{DOF})^2+(G_{DOF})^2\times(2\times 6 -1)+6\times G_{DOF}\times(2\times 6-1))$ FLOPs.  Finally, the cost of evaluating the stiffness matrix is given by the sum of FLOPs in the execution of Eqs.~\eqref{eq:comp_cost1} and \eqref{eq:stiff}: $N_e\times N_{GP}\times (G_{DOF})^2+(G_{DOF})^2\times(2\times 6 -1)+6\times G_{DOF}\times (2\times 6-1)+ N_{GJP}\times (2 \times 6 \times G_{DOF}))$. 
Assuming that the mesh is sufficiently refined, the overall computational complexity of this step is given to be $\mathcal{O}(G_{DOF}^3)$. The FLOPs for Eqs.~\eqref{eq:comp_cost1} and \eqref{eq:stiff} are $\mathcal{O}(G_{DOF}^2)$ and $\mathcal{O}(G_{DOF}^3)$, respectively. For a sufficiently refined distribution of nodes, the latter significantly dominates over the former ($\mathcal{O}(G_{DOF}^3) > \mathcal{O}(G_{DOF}^2)$) to give the cost of 2D f-EFG to be $\mathcal{O}(G_{DOF}^3)$.

It should be noted that the f-EFG method incurs a lower computational cost compared to the traditional FE solvers used for FDEs. In FE formulations, nonlocal analysis incurs an overall computational cost of $\mathcal{O}((\mathcal{N}_x+\mathcal{N}_y)\times(\bar{G}_{DOF})^3)$, where $(\mathcal{N}_x+\mathcal{N}_y)$ represents the number of elements within the nonlocal horizon at each evaluation point. This contrasts to the computational cost of $\mathcal{O}(G_{DOF}^3)$ noted here for f-EFG. Recall the definitions for $\mathcal{N}_{x}$ and $\mathcal{N}_{y}$ in Eq.~\eqref{eq: dynamic_rate}, and the choice of $\mathcal{N}_{x}=\mathcal{N}_{y}=12$ identified from the convergence studies. Furthermore, the number of degrees of freedom for the FE can be higher ($\bar{G}_{DOF}=6\times n$; $n$ being the number of FE nodes) compared to the EFG (${G}_{DOF}=3\times n$; $n$ being the number of EFG nodes). This is due to additional variables introduced in interpolation via FE to achieve $C^1$ continuity. Such an increase in continuity via EFG is achieved trivially without any increase in the number of degrees of freedom. Note that the cost is proportional $G_{DOF}^3$. Therefore, $\bar{G}_{DOF}=2G_{DOF}$ results in an order difference in overall cost between FE and EFG.   

Clearly, the FLOPs of the meshfree f-EFG solver is significantly lower ($<5\%$) compared to the mesh-based f-FE solver developed in \cite{patnaik2021fractional}. Importantly, the FLOPs of the f-EFG do not depend on the size of the nonlocal influence horizon. This is in contrast to FEM that demands highly refined meshes to achieve accurate results and therefore offers a significant advantage for small nonlocal length scales $h_l$. 

\section{Conclusion}
This work presents a meshfree computational method for solving 2D fractional-order partial differential equations (FDEs). The f-EFG method, integrated with Moving Least Squares (MLS) interpolants, previously developed for 1D problems is extended here for 2D FDEs. Thus, a 2D f-EFG solver is developed here for FDEs. Owing to its meshfree nature and global interpolation capabilities, the proposed method is particularly well-suited for handling spatial fractional derivatives, which inherently involve nonlocal differ-integral formulations across the domain. To evaluate the performance of the 2D f-EFG approach, we consider benchmark problems of nonlocal Kirchhoff plates modeled using fractional-order constitutive relations. The 2D f-EFG solver is verified for convergence and validated with benchmark results in the literature. Owing to the flexibility of the meshfree method, the 2D f-EFG solver demonstrates efficacy for non-uniform node distributions. Further, 2D f-EFG solver developed here presents with numerical solution of FDEs over irregular geometries, which are otherwise extremely challenging to obtain via traditional mesh-based approaches. Finally, the 2D f-EFG method demonstrates superior computational efficiency compared with 2D FE solvers for FDEs. 
In conclusion, the 2D f-EFG method offers a robust, accurate and computationally efficient alternative to mesh-based techniques for a numerical solution of FDEs. Its general formulation allows for easy extension to complex geometries and multiphysics interactions, making it a powerful and flexible framework for the numerical analysis of advanced materials and structures.\\

\noindent
\textbf{Acknowledgment}: S.D. and S.S. acknowledge the financial support from the Anusandhan National Research Foundation (ANRF), erstwhile Science and Engineering Research Board (SERB), India, under the startup research grant program (SRG/2022/000566). M.H.R. and S.S. acknowledge the financial support provided by the Defence Research and Development Organization (DRDO), India under the grant vide file number DMRL/DMR-326/A/TC.

\bibliographystyle{naturemag}
\bibliography{references}
\end{document}